\documentclass[12pt]{amsart}
\usepackage{epic,eepic}
 \newlength{\picunit}               
 \newlength{\templengthhoriz}        
 \newlength{\templengthvert}         
 \newlength{\temprule}               
 \newcount{\numlines}                
\newcommand{\encapsfig}[3]{                      
    \setbox0=\hbox{
	   \setlength{\unitlength}{#2\picunit}
	    \input #1.tex }              
   
     \templengthvert=\ht0                        
     \advance \templengthvert by \dp0            
     \advance \templengthvert by  10\unitlength  
     \advance \templengthvert by 0.5ex           
     \templengthhoriz=\wd0                       
     \advance \templengthhoriz by 0.75em         
     \kern 10\unitlength                         
     \if#3l                                      
     \box0                                       
     \kern -\templengthvert                      
     \numlines=0                               
     \loop \ifdim0pt<\templengthvert \advance\numlines by 1
				  \advance\templengthvert by -\baselineskip
				  \repeat
     \hangindent=\templengthhoriz                
     \hangafter=-\numlines                       
     \else                                       
     \temprule=\textwidth                        
     \advance \temprule by -\wd0                 
     \noindent\kern\temprule\box0                

     \kern -\templengthvert                      
     \numlines=0
     \loop \ifdim0pt<\templengthvert \advance\numlines by 1
				  \advance\templengthvert by -\baselineskip
				  \repeat
     \hangindent=-\templengthhoriz               
     \hangafter=-\numlines                       
     \fi                                         
   }

\newcommand{\point}{\vspace{3mm}\par\refstepcounter{subsection}{\bf \thesubsection.} }
\newcommand{\tpoint}[1]{\vspace{3mm}\par\refstepcounter{subsection}{\bf \thesubsection.} 
  {\em #1. ---} }
\newcommand{\epoint}[1]{\vspace{3mm}\par\refstepcounter{subsection}{\bf \thesubsection.} 
  {\em #1.} }
\newcommand{\bpoint}[1]{\vspace{3mm}\par\refstepcounter{subsection}{\bf \thesubsection.} 
  {\bf #1.} }

\newlength{\baseunit}               
\newcommand{\bpf}{\noindent {\em Proof.  }}
\newcommand{\epf}{\qed \vspace{+10pt}}

\newcommand{\zed}{\mathbb{Z}}

\newcommand{\com}{\mathbb{C}}

\newcommand{\De}{\Delta}

\newcommand{\proj}{\mathbb P}
\newcommand{\aff}{\mathbb A}

\newcommand{\eff}{\mathbb F}

\newcommand{\oh}{{\mathcal{O}}}
\newcommand{\cC}{{\mathcal{C}}}

\newcommand{\cK}{{\mathcal{K}}}
\newcommand{\cL}{{\mathcal{L}}}

\newcommand{\cU}{{\mathcal{U}}}

\newcommand{\cm}{{\mathcal{M}}}

\newcommand{\tC}{\tilde{C}}

\newcommand{\al}{\alpha}
\newcommand{\om}{\omega}
\newcommand{\be}{\beta}
\newcommand{\ga}{\gamma}

\newcommand{\de}{\delta}

\newcommand{\degen}{\operatorname{degen}}
\newcommand{\idim}{\operatorname{idim}}

\newcommand{\Sym}{\operatorname{Sym}}
\newcommand{\mbar}{\overline{M}}
\newcommand{\cmbar}{\overline{\cm}}
\newcommand{\mgt}{\cmbar_g(\proj^2,d)^+}

\newcommand{\mts}{\cmbar_3(\proj^2,4)^*}
\newcommand{\mtp}{\cmbar_3(\proj^2,4)^+}
\newcommand{\mos}{\cmbar_1(\proj^2,3)^*}
 
\begin{document}
\pagestyle{plain}
\title{The characteristic numbers of quartic plane curves}
\author{Ravi Vakil}
\date{October 23,  1998.}

\begin{abstract}
The characteristic numbers of smooth plane quartics are computed using
intersection theory on a component of the moduli space of stable maps.
This completes the verification of Zeuthen's prediction of
characteristic numbers of smooth plane curves.  A short sketch of a
computation of the characteristic numbers of plane cubics is also
given as an illustration.  
\end{abstract}
\maketitle

\tableofcontents

\section{Introduction}

\begin{verse}
The nineteenth century work on finding the characteristic numbers of
families of curves of higher degree is rich and lovely.  Understanding
it well enough to vindicate it and continue it, is possibly the most
important part of Hilbert's fifteenth problem remaining open. --- \cite{fkm} p. 193 
\end{verse}

\point
The classical {\em characteristic number problem} for smooth plane
curves (studied by Chasles, Zeuthen, Schubert, and other geometers of
the nineteenth century) is: {how many smooth degree $d$ plane
curves curves are there through $a$ fixed general points, and tangent
to $b$ fixed general lines (if $a+b=\binom {d+2} 2$)?} The success of earlier
geometers at correctly computing such numbers (and others from
enumerative geometry), despite the lack of a  firm theoretical foundation, led
Hilbert to include the justification of these methods as one of his
famous problems.  (For a more complete introduction to the history of
such problems, see \cite{k} and S.  Kleiman's introduction to
\cite{schubert}.)

H.G. Zeuthen predicted the characteristic numbers of curves of degree
at most 4.  Only with the advent of Fulton-Macpherson intersection
theory have these numbers begun to be verified.  The characteristic
numbers of the complete cubics were rigorously calculated by P. Aluffi
(\cite{acubics}) and S. Kleiman and R. Speiser (\cite{ksp}), and the
first ten characteristic numbers of the smooth quartics were computed
by Aluffi (\cite{aquartics}) and van Gastel (\cite{vangastel}).

It is interesting to compare the results and calculations with those
of Zeuthen (Section \ref{zsection}).  Although we use quite a different compactification,
unlike many other modern solutions of classical enumerative geometry
problems (such as the charactersitic numbers for twisted cubics),
the calculations are "similar".  They give a modern verification,
not only of these classical numbers, but, at least to some extent,
also of a classical approach. 

\epoint{Sketch of method}

The classical approach is to interpret the problem as the intersection
of divisors (corresponding to the incidence and tangency conditions) on
the parameter space of smooth curves, an open subvariety of a
projective space.  A ``good'' compactification must be given
(hopefully smooth, e.g. \cite{acubics}, at least in codimension 1),
and it must be checked that there are natural divisors on the
compactification that intersect (transversely) only in the open set
corresponding to smooth curves.

The method used here is as follows.  Kontsevich's moduli space of
stable maps gives us a compactification of the space of smooth
quartics.  (Explicitly: take the normalization of the component of the
moduli stack corresponding generically to closed immersions of smooth
curves.)  This compactification is birational to the parameter space
$\proj^{14}$ of smooth quartics, on which we have divisors $\al'$
(corresponding to curves through a fixed point), $\be'$ (corresponding
to curves tangent to a fixed line), and $\De$ (corresponding to nodal
curves), and 
\begin{equation}
\label{introeq}
\be' = 6 \al', \; \; \De = 27 \al'.
\end{equation} 

There are analogous divisors $\al$, $\be$, $\De_0$ on the
compactification, and equations (\ref{introeq}) remain true when
``lifted'' to the compactification, modulo ``boundary divisors''.  The
relevant boundary divisors are determined (Sections
\ref{boundarydivisors} and \ref{proofofugly}), and many of their
co-efficients in the ``lifts'' of (\ref{introeq}) are found using
one-parameter test families (Section \ref{testfamilies}).  The
intersections of the boundary divisors with cycles of the form $\al^a
\be^{13-a}$ ($0 \leq a \leq 13$) are calculated (up to two unknowns,
Section 7).  Then (the ``lifts'' of) the equations (\ref{introeq}) are
intersected with $\al^a \be^{13-a}$, giving a large number of linear
equations in the unknowns (including the characteristic numbers),
which can be solved (Section
\ref{linearalgebra}).  The characteristic numbers agree with Zeuthen's
predictions.  For example, there are 23,011,191,144 smooth plane
quartics tangent to 14 general lines.  

Section \ref{cubicaside} is a self-contained example of this approach,
giving a sketch of a quick calculation of the characteristic numbers
of smooth plane cubics.

\point
In summary, this paper resolves a problem of long-standing interest by
a classical approach, but using beautiful modern machinery, the theory
of stable maps.  If the measure of a new idea is its ability to shed
light on areas of previous interest, then this is yet another example
of the power of Kontsevich's moduli space of stable maps.

\epoint{Acknowledgements}  The author is grateful to W. Fulton and
B. Hassett for independently suggesting this problem, and for useful
comments.  He also thanks
A. J. de Jong and T. Graber for many fruitful discussions, and J.
Black and V. Vourkoutiotis for translating parts of \cite{schubert}.
Fulton and J. Harris made a helpful suggestion that substantially
changed the presentation of the argument.  The argument and ideas
presented owe much to the work (both published and unpublished) of P.
Aluffi, and his assistance throughout this project (including
providing a copy of \cite{zeuthen}) has been invaluable.

\section{Conventions and background results}

\point 
We follow the same conventions as in \cite{ell}.  We work over a fixed
algebraically closed field $k$ of characteristic 0.  By {\em scheme}, we
mean scheme of finite type over $k$.  By {\em variety}, we mean a
separated integral scheme.  By {\em stack} we mean Deligne-Mumford stack of
finite type over $k$.  All morphisms of schemes and stacks are assumed to be
defined over $k$, and fibre products are over $k$ unless otherwise
specified.

If $S$ is a Deligne-Mumford stack, then a {\em family of nodal curves}
(or a {\em nodal curve}) over $S$ is defined as usual (see \cite{ell}
2.2 for example; see \cite{dm} for the canonical treatment).

For basic definitions and results about {\em maps of nodal curves} and {\em stable maps}, see
\cite{fp} (or the brief summary in \cite{ell} 2.6).    
Let
${\cmbar_g(\proj^2,d)}$ be the stack whose category of sections of a
scheme $S$ is the category of families of stable maps to $\proj^2$
over $S$ of degree $d$ and arithmetic genus $g$.  For definitions and
basic results, see \cite{fp}.  It is a fine moduli stack of
Deligne-Mumford type.  There is a ``universal map'' over
${\cmbar_g(\proj^2,d)}$ that is a family of maps of nodal curves.
There is an open substack $\cm_g(\proj^2,d)$ that is a fine moduli
stack of maps of {\em smooth} curves.  There is a unique component of ${\cmbar_g(\proj^2,d)}$
that is the closure of such maps (of dimension $3d+g-1$);
call this component $\mgt$.  

If $\rho: C \rightarrow \proj^2 \times S$ is a family of maps of nodal
curves to $\proj^2$ over $S$, where $S$ is a Deligne-Mumford stack of
pure dimension $d$, then two classes $\al$ and $\be$ in $A^1 S$ (the
operational Chow ring of $S$), functorial in $S$, were defined in
\cite{ell} Section 3.  The divisor $\al$ corresponds to maps through a fixed general point, and
$\be$ corresponds to maps tangent to a fixed general line.  We say
that $\al^a \be^b[S]$ ($a+b=d$) are the {\em characteristic numbers}
of the family of maps.  If all the characteristic numbers of the
family are 0, we say the family is {\em enumeratively irrelevant}.  
Recall conditions (*) and (**) on families of maps of nodal curves, from
\cite{ell} Section 2.4:
\begin{itemize}
\item[(*)] Over a dense open substack of $S$, the curve $C$ is smooth, and $\rho$ factors 
$C \stackrel {\al} {\rightarrow} C' \stackrel {\rho'} {\rightarrow} \proj^2 \times S$ 
where $\rho'$ is unramified and gives a birational
map from $C'$ to its image; $\al$ is a degree $d_{\al}$ map with only simple ramification (i.e.
reduced ramification divisor); and 
the images of the simple
ramifications are distinct in $\proj^2$. 
\item[(**)]  Over the normal locus (a dense open substack) of $S$, each component of the normalization of $C$ (which is a family of maps
of nodal curves) satisfies (*).
\end{itemize}
If the family satisfies condition (**), then the characteristic numbers can be
interpreted enumeratively using \cite{ell} Theorem 3.15, as counting
maps (with multiplicity).

The classical {\em characteristic number problem} for curves in
$\proj^2$ studied by geometers of the last century is: {\em how
many irreducible nodal degree $d$ geometric genus $g$ maps are there
through $a$ general points, and tangent to $b$ general lines (if
$a+b=3d+g-1$)?} By \cite{ell} Theorem 3.15, this number is the degree of $\al^a
\be^b [\mgt]$.

\epoint{Genus 3 curves}
Let $\mts$ be the normalization of $\mtp$.

On the Deligne-Mumford stack $\cmbar_3$, let $h$ be the divisor that is the class of (the
closure of the locus corresponding to) smooth hyperelliptic curves.
Let $\de_0$ the the divisor corresponding to irreducible nodal curves.
Let $\de_1$ be the divisor corresponding to nodal curves with a component of arithmetic genus 1.

\section{Aside:  The complete cubics revisited}
\label{cubicaside}

\point As an example of the method, we sketch a derivation of the
characteristic numbers of smooth plane cubics.  The characteristic
numbers of smooth plane cubics were predicted by Zeuthen in the last
century.  They have since been calculated rigorously in the 1980's by
Aluffi (\cite{acubics}, using a smooth compactification, the complete
cubics) and Kleiman and Speiser (\cite{ksp}, using codimension 1
degenerations), and more recently by the author (\cite{ell}) and
Graber and Pandharipande (using the theory of gravitational
descendants,
\cite{gp}).  The numbers have also been computed (although not
rigorously proved) by degeneration of the point and tangency
conditions.  

\point Many
verifications will be left to the reader.  As an exercise, the reader
may enjoy using the same method to quickly calculate characteristic
numbers of smooth plane conics.  (In this case, the method turns out
to be identical in substance to the method of complete conics.)

Let $\mos$ be the normalization of the component of the moduli stack
$\cmbar_1(\proj^2,3)$ that is the closure of the locus of immersions
of smooth curves.  Then there are three enumeratively relevant
boundary divisors: 
\begin{enumerate}
\item[(i)] $\De_0$ is the closure of the locus of immersions
of nodal cubics,
\item[(ii)] $I$ is the closure of the locus of 3-to-1 maps from 
a smooth elliptic curve onto a line in $\proj^2$ (ramifying at 6 points), and 
\item[(iii)] $T$ is the closure of the locus of maps
from curves $C_0 \cup C_1$ where $C_i$ is smooth of genus $i$, the two curves meet at a node,
$C_0$ maps to a line, and $C_1$ maps 2-to-1 onto a line.
\end{enumerate}
(There are three other, enumeratively irrelevant, divisors; see
\cite{ratell} Lemma 3.14 for their description.)  The divisor $\De_0$
won't concern us in this example, but its analogue will be necessary
for the quartic case.

As described in Section 2 and \cite{ell} Section 3.16, there are also two divisors
$\al$ and $\be$ such that the characteristic number of cubics through
$a$ points and tangent to $b$ lines is the degree of $\al^a \be^b
[\mos]$.  

On the $\proj^9$ parametrizing plane cubics, there are analogous
divisors $\al'$ and $\be'$, and $\be' = 4 \al'$.  The two spaces are
birational, with isomorphic open subschemes parametrizing closed immersions.  Hence,
in $A_8(\mos)$, modulo enumeratively irrelevant divisors, 
\begin{equation}
\label{cubictake1}
4 \al = \be + t T + i I.
\end{equation}
for some rational $t$ and $i$.  We can find $t$ and $i$ by intersecting this with suitable
one-parameter families.  

Consider a pencil joining a general cubic curve and a triple line.  In
other words, if $p(x,y,z)=0$ describes a general cubic, consider the pencil
$\lambda p(x,y,z) + \mu x^3 =0$ with $[\lambda,\mu] \in \proj^1$. This
describes a family of nodal curves except at the point corresponding to
the triple line; perform stable reduction (for maps) to complete the family,
and get a map $\proj^1 \rightarrow \mos$.  On this family, compute
that $T=0$ (as the family misses $T$), $\al=1$ (because it's a pencil),
and $\be = 2$.  Check that the family intersects the Weil divisor $I$
transversely at one point, and hence $I=1/3$.  (It is essential to work with
stacks rather than schemes!  The 1/3 comes from the fact that the 
limit stable map has an automorphism group of order 3.)  Hence $i=6$.

Next, take a pencil joining a general cubic and $x^2 y = 0$.    On this 
family, $I=0$, $\al=1$, $\be = 3$, and $T = 1/2$.  Hence $t=2$, and (\ref{cubictake1}) can be rewritten 
\begin{equation}
\label{cubic}
4 \al = \be + 2 T + 6 I.
\end{equation}

\point \label{triplecover}
We can easily compute the characteristic numbers of $T$ and $I$.  For
example, the degree of $\be^8 [I]$ (the number of maps in $I$ tangent
to 8 fixed general lines) can be computed as follows.  For a map in
$I$ to be tangent to 8 general lines, the image line $\ell$ of the map
must pass through the intersections of two pairs of these lines (there
are $3 \binom 8 4 = 210$ choices of two pairs).  Then the map must be
a triple cover of $\ell$, branched over the intersection of $\ell$
with the 8 lines (which are 6 points).  The number of connected triple
covers of $\ell$ with 6 given branch points is $(\frac {3^5-3} {3!} )
= 40$.  (Proof: Rigidify the combinatorial problem by fixing some
other point in $\proj^1$, and labeling the 3 points mapping to it.
Monodromy about the six branch points gives transpositions in $S_3$,
and the product of these transpositions must be the identity.
Conversely, six such transpositions uniquely determine a cover, by the
Riemann existence theorem.  Thus five of the transpositions can be
chosen arbitrarily ($3^5$ choices), and the sixth is then determined.
However, the five cannot all be the same transposition, as then the
cover would be disconnected (leaving $3^5-3$ choices).  Finally,
divide by $3!$ to account for the labeling of the 3 points.)  Hence
the degree of $\be^8 [I] = 210 \times 40 = 8400$.  
(This is actually a special case of a formula of Hurwitz, \cite{hurwitz}.)  
The other non-zero
characteristic numbers of $I$ are $\al^2 \be^6[I] = 360$ and $\al
\be^7[I] = 2520$.

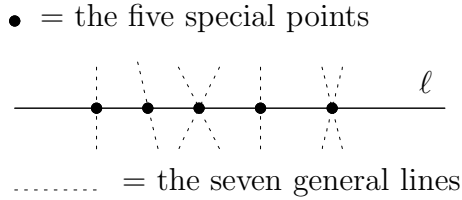
\begin{figure}
\begin{center}

	   \setlength{\unitlength}{.1\picunit}
	    \begingroup\makeatletter\ifx\SetFigFont\undefined%
\gdef\SetFigFont#1#2#3#4#5{%
  \reset@font\fontsize{#1}{#2pt}%
  \fontfamily{#3}\fontseries{#4}\fontshape{#5}%
  \selectfont}%
\fi\endgroup%
{\renewcommand{\dashlinestretch}{30}
\begin{picture}(6423,2609)(0,-10)
\put(1283,1236){\blacken\ellipse{150}{150}}
\put(1283,1236){\ellipse{150}{150}}
\dashline{60.000}(1283,1836)(1283,636)
\put(2033,1236){\blacken\ellipse{150}{150}}
\put(2033,1236){\ellipse{150}{150}}
\put(2783,1236){\blacken\ellipse{150}{150}}
\put(2783,1236){\ellipse{150}{150}}
\put(3683,1236){\blacken\ellipse{150}{150}}
\put(3683,1236){\ellipse{150}{150}}
\put(4733,1236){\blacken\ellipse{150}{150}}
\put(4733,1236){\ellipse{150}{150}}
\put(83,2511){\blacken\ellipse{150}{150}}
\put(83,2511){\ellipse{150}{150}}
\dashline{60.000}(1883,1911)(2183,636)
\dashline{60.000}(3083,1836)(2483,636)
\dashline{60.000}(2483,1836)(3083,636)
\dashline{60.000}(3683,1836)(3683,636)
\path(83,1236)(6383,1236)
\path(83,1236)(6383,1236)
\dashline{60.000}(4883,1836)(4583,636)
\dashline{60.000}(4583,1836)(4883,636)
\dashline{60.000}(83,36)(1283,36)
\put(458,2436){\makebox(0,0)[lb]{\smash{{{\SetFigFont{12}{14.4}{\rmdefault}{\mddefault}{\updefault}= the five special points}}}}}
\put(1658,36){\makebox(0,0)[lb]{\smash{{{\SetFigFont{12}{14.4}{\rmdefault}{\mddefault}{\updefault}= the seven general lines}}}}}
\put(6008,1461){\makebox(0,0)[lb]{\smash{{{\SetFigFont{12}{14.4}{\rmdefault}{\mddefault}{\updefault}$\ell$}}}}}
\end{picture}
} 
\end{center}
\caption{The five special points on $\ell$}
\label{cubicfig}
\end{figure}

\point
\label{cubict}
As another example, we compute the degree of $\al \be^7 [T]$ (the
number of maps in $T$ through a fixed point and tangent to 7 general
lines, with appropriate multiplicity).  For such maps, $C_0$ must go
through the fixed point.  Let $\ell$ be the image of $C_1$, and $m$ be
the image of $C_0$.  Then $\ell$ must pass through the intersection of
two pairs of the seven lines, so $\ell$ meets the seven lines at a
total of five ``special'' points.  (See Figure \ref{cubicfig} for a
pictorial representation.)  Through each of three of special points
one of the seven lines passes.  Through each of the other two special
points, two of the seven lines pass.

The cover $C_1 \rightarrow \ell$ ramifies at 4 of these 5 special points, and
$m$ meets $\ell$ at the fifth.  Each such map is counted with
multiplicity $2^a$, where $a$ is the number of the 7 lines passing
through the intersection of $m$ and $\ell$ (i.e. the image of the
node of the source curve).  There are $\binom 7 {2,2,3} / 2 = 105$
ways of choosing $\ell$.  Then $m$ can pass through one of the three
special points through which one of the seven lines pass (there are 3
ways of choosing this point, and the multiplicity is $2^1$), or $m$
can pass through one of the two special points through which two of
the seven lines pass (there are 2 ways of choosing this point, and the
multiplicity is $2^2$).

Hence the degree of $\al \be^7 [T]$ is $105(3 \times 2 + 2 \times 4) =
1470.$ (The other non-zero characteristic numbers of $T$ are $\al^4
\be^4 [T] = 24$, $\al^3 \be^5[T] = 240$, $\al^2 \be^6[T] = 885$.)

Thus the characteristic numbers of $I$ and $T$ can really be computed
by hand.

\point If we intersect (\ref{cubic}) with $\al^a \be^{8-a}$ ($0 \leq a \leq
8$), we have an equation relating two ``adjacent'' characteristic
numbers of smooth cubics, and characteristic numbers of $T$ and $I$.
As the degree of $\al^9 [\mos]$ is 1 (there is one smooth cubic through
9 general points), we can compute all the characteristic numbers
inductively.

As an example, the degree of $\al \be^8 [\mos]$ is 21004; from this we will calculate
the degree of $\be^9 [\mos]$.  Intersecting (\ref{cubic}) with $\be^8$, we get
$$
\deg( \be^9 [\mos]) = 4 \times 21004 - 2 \deg \be^8 [T] - 6 \deg \be^8 [I] .
$$
As $\be^8 [T] = 0$ (exercise) and $\deg \be^8 [I] = 8400$ from above,
$$
\deg( \be^9 [\mos]) = 4 \times 21004 - 6 \times 8400 = 33616.
$$

Thus the characteristic numbers of plane cubics can be really be computed by hand.
The moduli space of stable maps, by providing an excellent compactification
of the space of smooth cubics, makes this classical problem much easier.

\section{Boundary divisors}
\label{boundarydivisors}

\point 
\label{boundarydef}
We next describe the divisors on $\mts$ that are pertinent to the
argument.  Recall that a divisor $B$ is {\em enumeratively irrelevant}
if $\al^a \be^b [B]=0$ for all $a+b=13$.  We will see that the
following families are the enumeratively relevant boundary divisors on
$\mts$.  Each locus is clearly irreducible of dimension 13.  It will
not be immediate that these loci lie on $\mts$, but that will follow
from the rest of the argument (Theorem \ref{weilequation}).

Let $\De_0$ be the closure of the points in $\mts$ parametrizing
immersions of nodal curves.  Let $H$ be the closure of points
parametrizing smooth hyperelliptic curves mapping canonically to the
plane (and hence two-to-one onto a conic).  Let $I$ be the closure of
points parametrizing smooth genus 3 curves mapping canonically to a
line in the plane (i.e. $\rho^* \oh_{\proj^2(1)} \cong
\cK_C$).

\begin{figure}
\begin{center}

	   \setlength{\unitlength}{.1\picunit}
	    \begingroup\makeatletter\ifx\SetFigFont\undefined%
\gdef\SetFigFont#1#2#3#4#5{%
  \reset@font\fontsize{#1}{#2pt}%
  \fontfamily{#3}\fontseries{#4}\fontshape{#5}%
  \selectfont}%
\fi\endgroup%
{\renewcommand{\dashlinestretch}{30}
\begin{picture}(12987,8379)(0,-10)
\put(5175.000,7083.000){\arc{1950.000}{5.1072}{7.4592}}
\path(3075,3333)(4575,933)
\path(3075,933)(4575,3333)
\path(8625,933)(8625,2733)
\path(8550,1983)(9525,3333)
\path(8700,2358)(7800,3333)
\path(12075,8283)(12075,5883)
\path(12975,7983)(11175,7983)
\path(1125,8283)(1125,5883)
\path(2025,7983)(225,7983)
\path(5850,8283)(5850,5883)
\put(3525,33){\makebox(0,0)[lb]{\smash{{{\SetFigFont{12}{14.4}{\rmdefault}{\mddefault}{\updefault}$X$}}}}}
\put(4050,2058){\makebox(0,0)[lb]{\smash{{{\SetFigFont{12}{14.4}{\rmdefault}{\mddefault}{\updefault}$x$}}}}}
\put(4650,3483){\makebox(0,0)[lb]{\smash{{{\SetFigFont{12}{14.4}{\rmdefault}{\mddefault}{\updefault}$C_1$}}}}}
\put(4575,1233){\makebox(0,0)[lb]{\smash{{{\SetFigFont{12}{14.4}{\rmdefault}{\mddefault}{\updefault}$C_2$}}}}}
\put(4800,783){\makebox(0,0)[lb]{\smash{{{\SetFigFont{12}{14.4}{\rmdefault}{\mddefault}{\updefault}(2,1)}}}}}
\put(8475,33){\makebox(0,0)[lb]{\smash{{{\SetFigFont{12}{14.4}{\rmdefault}{\mddefault}{\updefault}$Y$}}}}}
\put(8100,2358){\makebox(0,0)[lb]{\smash{{{\SetFigFont{12}{14.4}{\rmdefault}{\mddefault}{\updefault}$y_2$}}}}}
\put(9375,2733){\makebox(0,0)[lb]{\smash{{{\SetFigFont{12}{14.4}{\rmdefault}{\mddefault}{\updefault}$C_3$ (1,0)}}}}}
\put(7275,3483){\makebox(0,0)[lb]{\smash{{{\SetFigFont{12}{14.4}{\rmdefault}{\mddefault}{\updefault}$C_2$ (1,0)}}}}}
\put(8700,1833){\makebox(0,0)[lb]{\smash{{{\SetFigFont{12}{14.4}{\rmdefault}{\mddefault}{\updefault}$y_3$}}}}}
\put(12000,4983){\makebox(0,0)[lb]{\smash{{{\SetFigFont{12}{14.4}{\rmdefault}{\mddefault}{\updefault}$Q$}}}}}
\put(12225,6333){\makebox(0,0)[lb]{\smash{{{\SetFigFont{12}{14.4}{\rmdefault}{\mddefault}{\updefault}$C_1$}}}}}
\put(12225,5958){\makebox(0,0)[lb]{\smash{{{\SetFigFont{12}{14.4}{\rmdefault}{\mddefault}{\updefault}(2,3)}}}}}
\put(10875,8208){\makebox(0,0)[lb]{\smash{{{\SetFigFont{12}{14.4}{\rmdefault}{\mddefault}{\updefault}$C_2$}}}}}
\put(10875,7608){\makebox(0,0)[lb]{\smash{{{\SetFigFont{12}{14.4}{\rmdefault}{\mddefault}{\updefault}(2,0)}}}}}
\put(975,4983){\makebox(0,0)[lb]{\smash{{{\SetFigFont{12}{14.4}{\rmdefault}{\mddefault}{\updefault}$T$}}}}}
\put(1275,6633){\makebox(0,0)[lb]{\smash{{{\SetFigFont{12}{14.4}{\rmdefault}{\mddefault}{\updefault}$C_1$}}}}}
\put(0,8133){\makebox(0,0)[lb]{\smash{{{\SetFigFont{12}{14.4}{\rmdefault}{\mddefault}{\updefault}$C_2$}}}}}
\put(1200,8058){\makebox(0,0)[lb]{\smash{{{\SetFigFont{12}{14.4}{\rmdefault}{\mddefault}{\updefault}$t$}}}}}
\put(0,7608){\makebox(0,0)[lb]{\smash{{{\SetFigFont{12}{14.4}{\rmdefault}{\mddefault}{\updefault}(1,0)}}}}}
\put(5775,4983){\makebox(0,0)[lb]{\smash{{{\SetFigFont{12}{14.4}{\rmdefault}{\mddefault}{\updefault}$P$}}}}}
\put(6000,8133){\makebox(0,0)[lb]{\smash{{{\SetFigFont{12}{14.4}{\rmdefault}{\mddefault}{\updefault}$C_1$ (2,2)}}}}}
\put(6225,6708){\makebox(0,0)[lb]{\smash{{{\SetFigFont{12}{14.4}{\rmdefault}{\mddefault}{\updefault}$C_2$ (2,0)}}}}}
\put(1350,6183){\makebox(0,0)[lb]{\smash{{{\SetFigFont{12}{14.4}{\rmdefault}{\mddefault}{\updefault}(3,3)}}}}}
\put(12150,8058){\makebox(0,0)[lb]{\smash{{{\SetFigFont{12}{14.4}{\rmdefault}{\mddefault}{\updefault}$q$}}}}}
\put(4725,3033){\makebox(0,0)[lb]{\smash{{{\SetFigFont{12}{14.4}{\rmdefault}{\mddefault}{\updefault}(2,2)}}}}}
\put(8850,933){\makebox(0,0)[lb]{\smash{{{\SetFigFont{12}{14.4}{\rmdefault}{\mddefault}{\updefault}$C_1$ (2,3)}}}}}
\end{picture}
} 
\end{center}
\caption{Source curves corresponding to general points of the boundary divisors $T$, $P$, $Q$, $X$, $Y$ (components labeled (degree, genus))}
\label{tpqxy}
\end{figure}
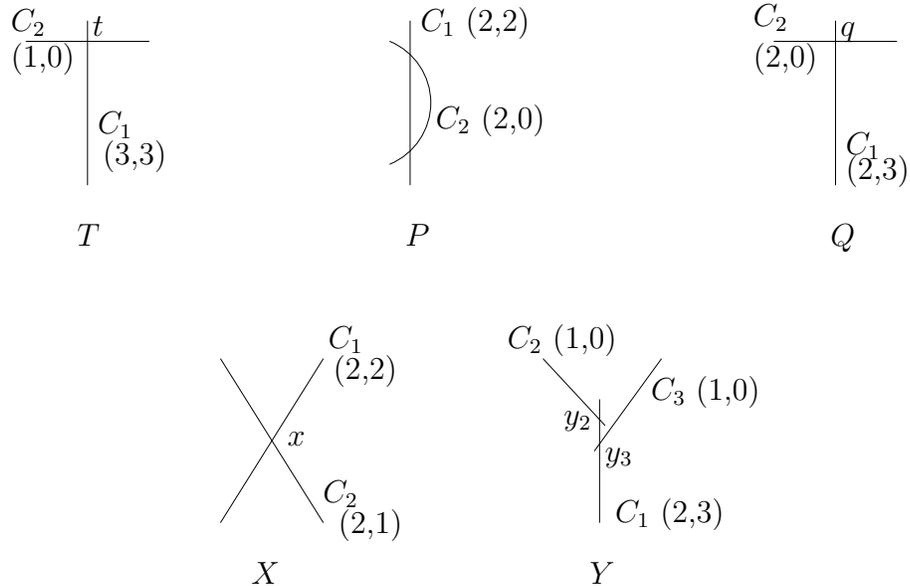

The boundary divisors $T$, $P$, $Q$, $X$, and $Y$ are described in
Figure \ref{tpqxy}.  The source curve is given (diagramatically),
where the components are labeled $C_i$, and each component is
labeled with an ordered pair of the degree and genus of the map
(restricted to that component).  For $T$, the component $C_1$ is
mapped to $\proj^2$ by the line bundle $\cK_{C_1}(-t)$.  (Equivalently, if
$\rho$ is the morphism from $C_1$ to $\proj^2$, $\rho^{-1}
\oh_{\proj^2}(1) \cong \cK_{C_1}(-t)$.)  The 
component $C_1$ triple-covers a line.  For $P$, $Q$, $X$, and $Y$, the
image of $C_1$ is necessarily a double-line.  For $Q$, the point $q$
is required to be a Weierstrass point of $C_1$, and the image of $C_2$
(a smooth plane conic) is required to be tangent to the image of
$C_1$.  For $X$, the map from $C_2$ to $\proj^2$ (a double cover of a
line) is required to ramify at the point $x$.  For $Y$, the points
$y_2$ and $y_3$ are required to be hyperelliptically conjugate.

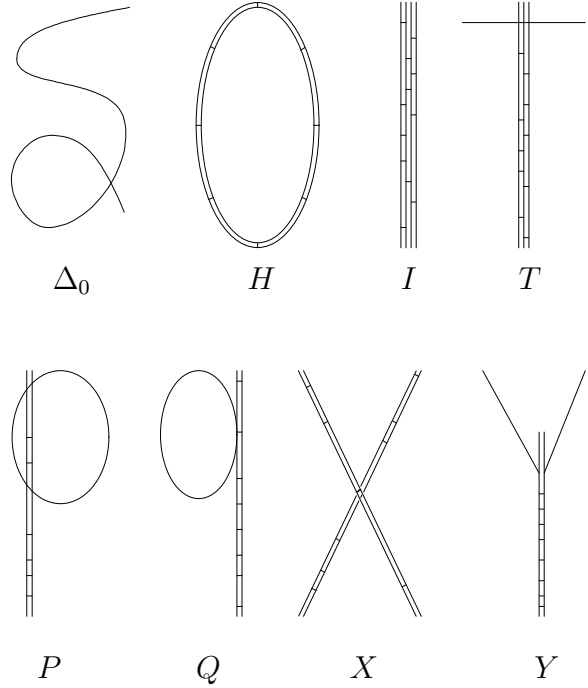
\begin{figure}
\begin{center}

	   \setlength{\unitlength}{.1\picunit}
	    \begingroup\makeatletter\ifx\SetFigFont\undefined%
\gdef\SetFigFont#1#2#3#4#5{%
  \reset@font\fontsize{#1}{#2pt}%
  \fontfamily{#3}\fontseries{#4}\fontshape{#5}%
  \selectfont}%
\fi\endgroup%
{\renewcommand{\dashlinestretch}{30}
\begin{picture}(8434,9963)(0,-10)
\put(3622,8136){\ellipse{1800}{3600}}
\put(3622,8136){\ellipse{1650}{3450}}
\path(5722,9936)(5722,6336)
\path(5797,9936)(5797,6336)
\path(5872,9936)(5872,6336)
\path(5947,9936)(5947,6336)
\path(7522,9936)(7522,6336)
\path(7597,9936)(7597,6336)
\path(7447,9936)(7447,6336)
\path(6622,9636)(8422,9636)
\put(734,3561){\ellipse{1424}{1950}}
\path(322,4536)(322,936)
\path(247,4536)(247,936)
\put(2759,3599){\ellipse{1126}{1874}}
\path(3322,936)(3322,4536)
\path(3397,4536)(3397,936)
\path(7822,936)(7822,3636)
\path(7747,936)(7747,3636)
\path(7747,3036)(6922,4536)
\path(7822,3036)(8422,4536)
\path(5722,9636)(5797,9636)
\path(5797,9111)(5872,9111)
\path(5872,9411)(5947,9411)
\path(5722,6636)(5797,6636)
\path(5947,7011)(5872,7011)
\path(5797,7311)(5872,7311)
\path(5722,7986)(5797,7986)
\path(5872,8886)(5947,8886)
\path(5722,7611)(5797,7611)
\path(5872,8661)(5797,8661)
\path(5722,8436)(5797,8436)
\path(5872,8286)(5947,8286)
\path(7447,9186)(7522,9186)
\path(7522,8886)(7597,8886)
\path(7522,8436)(7597,8436)
\path(7447,8211)(7522,8211)
\path(7522,6486)(7597,6486)
\path(7447,6786)(7522,6786)
\path(7597,7236)(7522,7236)
\path(7447,7461)(7522,7461)
\path(7447,7761)(7522,7761)
\path(7522,7986)(7597,7986)
\path(247,1236)(322,1236)
\path(247,1761)(322,1761)
\path(247,2136)(322,2136)
\path(247,3561)(322,3561)
\path(247,3186)(322,3186)
\path(247,1536)(322,1536)
\path(3322,3636)(3397,3636)
\path(3322,1086)(3397,1086)
\path(3322,1536)(3397,1536)
\path(3322,1836)(3397,1836)
\path(3322,2211)(3397,2211)
\path(3322,2586)(3397,2586)
\path(3322,2961)(3397,2961)
\path(3322,4386)(3397,4386)
\path(7747,1086)(7822,1086)
\path(7747,1236)(7822,1236)
\path(7747,1461)(7822,1461)
\path(7747,1761)(7822,1761)
\path(7747,2061)(7822,2061)
\path(7747,2286)(7822,2286)
\path(7747,2511)(7822,2511)
\path(7747,2736)(7822,2736)
\path(4402,4326)(4338,4292)
\path(4653,3793)(4593,3763)
\path(5797,1401)(5741,1371)
\path(5595,3808)(5658,3771)
\path(5460,3512)(5520,3478)
\path(4417,1352)(4481,1315)
\path(4548,1615)(4608,1581)
\path(4320,9272)(4256,9235)
\path(2932,9287)(2996,9257)
\path(2906,7048)(2977,7075)
\path(4338,7048)(4271,7078)
\path(4522,8136)(4447,8136)
\path(2722,8136)(2797,8136)
\path(3622,9936)(3622,9861)
\path(3622,6411)(3622,6336)
\path(4774,2079)(4828,2049)
\path(5069,2765)(5133,2791)
\path(5946,4532)(5137,2845)
\path(5132,2677)(5149,2714)
\path(5095,2755)(5115,2799)
\path(5073,2709)(4226,943)
\path(4222,4536)(5947,936)
\path(6022,936)(4297,4536)
\path(5928,4483)(5985,4453)
\path(4297,936)(5108,2630)
\path(5174,2765)(6022,4536)
\path(1747,9861)	(1690.321,9850.602)
	(1634.818,9840.251)
	(1580.482,9829.942)
	(1527.303,9819.672)
	(1475.271,9809.435)
	(1424.378,9799.228)
	(1374.614,9789.047)
	(1325.968,9778.886)
	(1278.433,9768.742)
	(1231.998,9758.611)
	(1186.655,9748.487)
	(1142.393,9738.368)
	(1099.203,9728.248)
	(1057.075,9718.123)
	(1016.001,9707.989)
	(975.971,9697.842)
	(936.975,9687.677)
	(899.004,9677.491)
	(826.099,9657.035)
	(757.181,9636.440)
	(692.174,9615.671)
	(631.002,9594.694)
	(573.589,9573.476)
	(519.862,9551.981)
	(469.742,9530.175)
	(423.157,9508.024)
	(380.029,9485.494)
	(340.283,9462.550)
	(303.844,9439.157)
	(240.583,9390.891)
	(189.643,9340.420)
	(150.419,9287.471)
	(122.306,9231.767)
	(104.701,9173.035)
	(97.000,9111.000)

\path(97,9111)	(99.263,9070.336)
	(109.653,9033.268)
	(152.771,8968.804)
	(222.278,8915.380)
	(265.653,8892.111)
	(314.095,8870.766)
	(367.097,8851.066)
	(424.147,8832.732)
	(484.736,8815.486)
	(548.354,8799.049)
	(614.493,8783.142)
	(682.641,8767.487)
	(752.289,8751.804)
	(822.929,8735.816)
	(894.049,8719.244)
	(965.141,8701.808)
	(1035.694,8683.230)
	(1105.200,8663.232)
	(1173.147,8641.534)
	(1239.028,8617.858)
	(1302.331,8591.926)
	(1362.548,8563.458)
	(1419.168,8532.177)
	(1471.682,8497.802)
	(1519.581,8460.056)
	(1562.354,8418.660)
	(1599.492,8373.335)
	(1630.486,8323.802)
	(1654.825,8269.784)
	(1672.000,8211.000)

\path(1672,8211)	(1683.958,8145.132)
	(1691.294,8078.385)
	(1694.165,8010.953)
	(1692.730,7943.032)
	(1687.148,7874.817)
	(1677.576,7806.504)
	(1664.173,7738.288)
	(1647.097,7670.364)
	(1626.506,7602.929)
	(1602.560,7536.177)
	(1575.415,7470.303)
	(1545.231,7405.504)
	(1512.165,7341.974)
	(1476.377,7279.909)
	(1438.024,7219.504)
	(1397.264,7160.955)
	(1354.256,7104.457)
	(1309.158,7050.206)
	(1262.129,6998.397)
	(1213.326,6949.225)
	(1162.908,6902.885)
	(1111.034,6859.574)
	(1057.862,6819.487)
	(1003.549,6782.818)
	(948.254,6749.764)
	(892.136,6720.519)
	(835.353,6695.280)
	(778.063,6674.241)
	(720.425,6657.599)
	(662.596,6645.547)
	(604.735,6638.282)
	(547.000,6636.000)

\path(547,6636)	(492.885,6641.079)
	(439.811,6655.289)
	(388.164,6677.868)
	(338.327,6708.054)
	(290.689,6745.083)
	(245.634,6788.194)
	(203.547,6836.623)
	(164.815,6889.609)
	(129.823,6946.388)
	(98.957,7006.197)
	(72.602,7068.275)
	(51.145,7131.859)
	(34.970,7196.186)
	(24.464,7260.494)
	(20.012,7324.019)
	(22.000,7386.000)

\path(22,7386)	(29.524,7442.454)
	(42.057,7498.519)
	(59.383,7553.738)
	(81.286,7607.651)
	(107.546,7659.801)
	(137.949,7709.729)
	(172.276,7756.978)
	(210.310,7801.087)
	(251.835,7841.601)
	(296.632,7878.060)
	(344.486,7910.006)
	(395.179,7936.980)
	(448.494,7958.525)
	(504.214,7974.182)
	(562.122,7983.493)
	(622.000,7986.000)

\path(622,7986)	(662.829,7984.477)
	(702.660,7981.299)
	(741.541,7976.415)
	(779.521,7969.773)
	(852.967,7951.010)
	(923.383,7924.596)
	(991.153,7890.121)
	(1056.663,7847.172)
	(1120.296,7795.338)
	(1182.438,7734.206)
	(1243.471,7663.365)
	(1273.693,7624.174)
	(1303.782,7582.402)
	(1333.787,7537.996)
	(1363.755,7490.905)
	(1393.734,7441.078)
	(1423.773,7388.463)
	(1453.920,7333.009)
	(1484.223,7274.664)
	(1514.729,7213.376)
	(1545.487,7149.095)
	(1576.545,7081.768)
	(1607.952,7011.344)
	(1639.754,6937.772)
	(1655.818,6899.789)
	(1672.000,6861.000)

\put(622,5736){\makebox(0,0)[lb]{\smash{{{\SetFigFont{12}{14.4}{\rmdefault}{\mddefault}{\updefault}$\Delta_0$}}}}}
\put(3472,5736){\makebox(0,0)[lb]{\smash{{{\SetFigFont{12}{14.4}{\rmdefault}{\mddefault}{\updefault}$H$}}}}}
\put(5722,5736){\makebox(0,0)[lb]{\smash{{{\SetFigFont{12}{14.4}{\rmdefault}{\mddefault}{\updefault}$I$}}}}}
\put(7447,5736){\makebox(0,0)[lb]{\smash{{{\SetFigFont{12}{14.4}{\rmdefault}{\mddefault}{\updefault}$T$}}}}}
\put(397,36){\makebox(0,0)[lb]{\smash{{{\SetFigFont{12}{14.4}{\rmdefault}{\mddefault}{\updefault}$P$}}}}}
\put(2722,36){\makebox(0,0)[lb]{\smash{{{\SetFigFont{12}{14.4}{\rmdefault}{\mddefault}{\updefault}$Q$}}}}}
\put(4972,36){\makebox(0,0)[lb]{\smash{{{\SetFigFont{12}{14.4}{\rmdefault}{\mddefault}{\updefault}$X$}}}}}
\put(7672,36){\makebox(0,0)[lb]{\smash{{{\SetFigFont{12}{14.4}{\rmdefault}{\mddefault}{\updefault}$Y$}}}}}
\end{picture}
} 
\end{center}
\caption{The images of the general maps in the boundary divisors $\De_0$, $H$, $I$, $T$, $P$, $Q$, $X$, $Y$} 
\label{prettypictures}
\end{figure}

Figure \ref{prettypictures} depicts images of the maps
corresponding to the general points of each of the divisors described
(with ramifications of the maps indicated suggestively).

The fundamental theorem of this section is the following.

\tpoint{Theorem}  {\em The enumeratively relevant boundary divisors of $\mts$ are $\De_0$, $H$, $I$, $T$, $P$, $Q$, $X$, $Y$.}
\label{ugly}

The proof is given in Section \ref{proofofugly}.   

\tpoint{Corollary}
{\em Modulo enumeratively irrelevant divisors, in $A_{13}(\mts)$, 
\begin{eqnarray}
6 \al &=& \be + hH + iI + tT + pP+qQ+ + xX + yY
\label{e1} \\
27 \al &=& \De_0+ h' H + i' I + t' T + p'P+q'Q+x' X + y'Y.
\label{e2}
\end{eqnarray}
for some rational numbers $h$, $i$, \dots, $x'$, $y'$.}

\bpf Let $\cU \subset \mts$  be the open subscheme corresponding to closed immersions of genus 3 curves.  Then
$\cU \cong \proj^{14} \setminus \Xi$ (where $\proj^{14}$ is the
Hilbert scheme parametrizing plane quartics, and $\Xi$ is a subset of
codimension greater than 1), as both sides represent the same functor.
Then by standard arguments, $\be|_{\cU} = 6 \al|_{\cU}$ and
$\De_0|_{\cU} = 27 \al|_{\cU}$ in $A_{13}(\cU)$, so in $A_{13}(\mts)$,
$\be = 6 \al$ and $\De_0 = 27 \al$ modulo boundary divisors except
$\De_0$. \epf

\epoint{A criterion for 1-parameter families to intersect divisors with multiplicity 1}
\label{node1}  
Let $\cC \rightarrow M$ be a family of nodal curves over a stack
$M$, such that the curve over the generic point is smooth.  Let $\De$
be an irreducible divisor on $M$ such that the universal curve over
$\De$ is singular (i.e. has a node).  Let $f: S \rightarrow M$ be a
morphism from a smooth curve to $M$, intersecting $\De$ at a
point $s \in S$, such that the pullback of the universal curve $\cC$
to the generic point of $S$ is smooth.  Suppose that $\De$ is locally Cartier at $f(s)$.  Recall that if the total space
of the pullback of the universal curve $\cC$ to $S$ is smooth above
$s$, then $f^* \De$ contains $s$ with multiplicity one, i.e. the
one-parameter family intersects $\De$ transversely. (Sketch of proof:
the formal deformation space of a node is smooth and one-dimensional;
let $(D,0)$ be this pointed space.  The universal curve over $D$ is
smooth, and the universal curve pulled back to a cover of $D$ ramified
at 0 is singular.  Choose any node of the curve above $s$.  Then the map
$S \rightarrow M$ induces a morphism $\pi$ from a formal neighborhood
of $s \in S$ to $D$, and as the total space of the universal curve
over $S$ is smooth above $s$, this map must be unramified, so $\pi$ is
\'{e}tale. But $f^{-1}(\De)$ is scheme-theoretically contained in
$\pi^* 0$ which is the reduced point $s$, so $f^{-1}(\De)$ is a
reduced point.)

\epoint{Description of $H$, $\De_0$, $X$ in terms of $h$, $\de_0$, $\de_1$}
\label{m3bar}
Let $\psi$ be the natural morphism $\mts \rightarrow \cmbar_3$.  Then
let $a$ be the multiplicity of $\psi^* h$ along $H$, $b$ be the
multiplicity of $\psi^* \de_0$ along $\De_0$, and $c$ be the
multiplicity of $\psi^* \de_1$ along $X$; $a$, $b$, and $c$ are
integers.  We will see later that $a=b=c=1$, using Criterion \ref{node1}.

\epoint{The enumerative geometry of $I$}
\label{tegoi}
A dimension count shows that of the 12-dimensional family of quadruple
covers of $\proj^1$ by smooth genus 3 curves, an 11-dimensional family
corresponds to canonical covers.  In other words, if 11 general points
are fixed on $\proj^1$, there are a finite number of quadruple
canonical covers branched at those 11 points; call this number
$\iota$.

Let $M$ be the space of genus 3 degree 4 admissible covers (with labeled
branch points), and $D_0$ the divisor that is the closure of the locus
of canonically mapped smooth curves.  If $\pi: M \rightarrow
\mbar_{0,12}$ is the natural map remembering only the branch points,
let $D = \pi_* D_0$.  (Surprisingly, $D$ has multiplicity 120; see
Section \ref{weirdIfact}.)  Let $\De_I$ be the boundary divisor on
$\mbar_{0,12}$ whose general point parametrizes a curve with two
components, with 2 of the marked points on one of the components, and
let $S_I$ be the set of boundary divisors not supported on $\De_I$.
Let $B$ be the one-parameter family described in the previous
paragraph (with 11 labeled points fixed and 1 moving), so $B
\cdot \De_I = 11$ and $B \cdot \De = 0$ for any $\De \in S$.  By
symmetry, $B$ meets each of the components of $\De_I$ with equal
multiplicity, and $D$ contains each of the components of $\De_I$ with
equal multiplicity.  Then as $B \cdot D =
\iota$, $D \equiv \frac \iota {11} \De_I \pmod S$.

\epoint{The enumerative geometry of $T$}
\label{tegot}
Similar to the previous case, consider $\mbar_{0,11}$, where the 11
points are labeled $u$, $p_1$, \dots, $p_{10}$.  Let $\De_T$ be the
boundary divisor where (generically) the curve has 2 components, one
with two points $p_i$, $p_j$ and one with the rest.  Let $\De_{T,u}$
be the boundary divisor where (generically) the curve has 2
components, one with two points $u$, $p_i$, and one with the rest.
Let $S$ be the set of boundary divisors not supported on $\De_T \cup
\De_{T,u}$.  Let $D$ be the divisor on $\mbar_{0,11}$ that is the
closure of the pushforward of the points of the pointed Hurwitz scheme
(where the marked points are a point $t$ and the branch points $p_1$,
\dots, $p_{10}$) corresponding to maps induced by the linear system
$\cK_C(-t)$.  

If 10 general points $p_i$ are fixed, then there are $(3^9-3)/3! =
3280$ possible connected triple covers branched there (see
Section \ref{triplecover} for an explanation of how to count connected triple
covers).  For each such cover $\pi: C
\rightarrow \proj^1$ there is exactly one point $t \in C$ such that
$\pi$ comes from the linear system $\cK_C(-t)$: if $|\cL|$ is the
linear system corresponding to $\pi$, then $h^0(C, \cL) \geq 2$, so by
Riemann-Roch, $h^0(C, \cK \otimes \cL^{-1}) \geq 1$.  But $h^0(C, \cK
\otimes \cL^{-1}) < 2$ as $\cK \otimes \cL^{-1}$ is a degree 1 line
bundle on an irrational curve, so $h^0(C, \cK \otimes \cL^{-1}) = 1$,
and $\cK \otimes \cL^{-1} \cong \oh(t)$ for a unique $t \in C$.

If 9 of the points $p_1$, \dots, $p_9$ and $u$ are fixed on $\proj^1$,
then let $\tau$ be the number of genus 3 triple covers $\pi: C
\rightarrow \proj^1$ branched at the 9 points $p_1$, \dots, $p_9$ (and
one other) with a point $t \in C$ with $\pi(t)=u$, such that $\pi$ is
induced by the linear series $\cK_C(-t)$.

If $B$ is the family described two paragraphs previously (with the
$p_i$ fixed and the $u$ moving) then $B \cdot \De_T = 0$,
$B \cdot \De_{T,u} = 10$, and $B \cdot D = 3280$.  If
$B'$ is the family described in the preceding paragraph (with $p_1$,
\dots, $p_9$ and $u$ fixed), 
then $B' \cdot \De_T = 9$, $B'
\cdot \De_{T,u} = 1$, and $B' \cdot D = \tau$.  Hence 
$D \equiv \left( \frac {\tau - 328} 9 \right)
\De_T + 328 \De_{T,u} \pmod S$.

\bpoint{Description of $I$ as a degeneracy locus}

Let $\pi: C \rightarrow S$
be a family of smooth genus 3 curves, and $\cL$ an invertible sheaf on
$C$ of (relative) degree 4, with sections $s_0$, $s_1$, $s_2 \in
h^0(C,\cL)$ giving a base-point free family of stable maps $C
\rightarrow \proj^2 \times S$.  This induces a morphism 
$S \rightarrow \mtp$.  Suppose this lifts to a morphism $\phi: S
\rightarrow \mts$ (e.g. if $S$ is normal), and suppose further that
$\phi(S)$ is not contained in $I$.  The subset of $S$ where the curve
maps to a line is a degeneracy locus (where the dimension of the
vector space spanned by $s_0$, $s_1$, and $s_2$ in a fiber is at most
2, \cite{f} Ch. 14).

\tpoint{Lemma} {\em 
If $m_{\degen}$ is the multiplicity with which an irreducible Weil
divisor $D$ appears in the degeneracy locus, and $m_I$ is the
multiplicity with which $D$ appears in $\phi^* I$, then $m_{\degen} =
m_I$.}
\label{iamdegenerate}

\bpf
If $S=\mts$ (with the universal family, and the sections given by
$s_0$, $s_1$, $s_2$ given by pullbacks of the co-ordinates $x$, $y$,
$z$ on $\proj^2$) and $D=I$, then $m_I = 1$, and $m_{\degen}$ is a positive
integer $k$.  By pulling back to an appropriate family, we see that
$k=1$ --- for example, fix a general genus 3 curve $C$ and 3 general
sections $s_0$, $s_1$, $s_2'$ of $\cK_C$, and consider the family $C
\times \aff^1 \rightarrow \proj^2 \times \aff^1$ (with co-ordinate $t$ on $\aff^1$) given by
$$
C \stackrel {(s_0, s_1, t s_2)} \longrightarrow \proj^2
$$
($t \in k$).  This family has $m_{\degen}  = m_I = 1$.

Finally, if $S$ is any other family of maps inducing a morphism $\phi: S
\rightarrow \mts$, then the degeneracy locus and $\phi^* I$ are both
pullbacks of the analogous loci on $\mts$, so $m_{\degen} = m_I$ on
this family as well.
\epf

\section{Proof of Theorem \ref{ugly}}
\label{proofofugly}

This proof is tedious and unenlightening, and the casual reader should
probably skip it.  

For simplicity, we divide the proof into a series
of steps.

\point  If $\cC \rightarrow \proj^2$ is a family of stable maps over $S$,
define the {\em intersection dimension} of the family (denoted $\idim S$) to be
the largest integer $n$ such that there is an integer $a$ ($0 \leq a
\leq n$) so there are maps in the family through $a$ fixed general
points and tangent to $n-a$ fixed general lines.  (Recall that a line
$\ell \subset \proj^2$ is tangent to a map $\rho: C \rightarrow
\proj^2$ if $\rho^* \ell$ is not a union of reduced points.)
Clearly $\idim(S) \leq \dim (S)$ (this is a consequence of
\cite{ell} Section 3; $\idim(S)$ is also bounded by the image of the
$S$ in the moduli space of stable maps).  Thus the Theorem asserts
that the only boundary divisors of $\mts$ that have intersection
dimension 13 are those listed.  For the rest of the proof, suppose
$\Xi$ is an irreducible boundary divisor of intersection dimension 13.

\point If $\cC \rightarrow \proj^2$ is a family of degree $d$ genus $g$
maps over an irreducible scheme $S$ ($1 \leq d \leq 4$, $0 \leq g \leq
3$) and the curve over a general $k$-point of $S$ is irreducible,
then it is easy to verify that the intersection dimension of the
family is at most that given in Table \ref{idimmax}, and that if
equality holds, then the generic source curve must be smooth.  Note that
if $d>g$ then the maximum is $3d+g-1$ (which is the virtual dimension
of the moduli space of degree $d$ genus $g$ stable maps to $\proj^2$).

\begin{table}
\begin{center}
\begin{tabular} {|r|c|c|c|c|}\hline
& $d=1$ & 2 & 3 & 4 \\ \hline
$g=0$ & 2 & 5 & 8 & 11 \\
1 & & 6 & 9 & 12 \\
2 & & 8 & 10 & 13 \\
3 & & 10 & 12 & 14  \\ \hline
\end{tabular}
\end{center}
\caption{Maximum intersection dimension of families of maps of irreducible curves}
\label{idimmax}
\end{table}

\point  \label{g3d4}
Suppose that the general (source) curve has a component of arithmetic
genus 3 that maps with degree 4.  Then this is the only component of
the general curve.  If the image of the general curve is reduced, then
(as the general map in $\Xi$ isn't an immersion of a smooth curve),
the image of $\Xi$ in $\proj^{14}$ must be the discriminant locus.
Then $\Xi=\De_0$.

If the image of the general curve is non-reduced, then it is either a
double conic or a quadruple line.  (As the general curve is
irreducible, in the first case the conic must be smooth.)  If the
general map is a double cover of a smooth conic, then $\Xi$ lies in
$H$.  As $\dim H = 13$ and $H$ is irreducible, $\Xi=H$.  If the
general map is a quadruple cover of a line, then as $\Xi \subset
\mts$, the general map is a limit of canonical maps.  As the general
curve in $\Xi$ is irreducible, the general map is given by the
canonical sheaf (i.e. the pullback of $\oh_{\proj^2}(1)$ is isomorphic to
the canonical sheaf), so $\Xi = I$ (as $\dim I = 13$ and $I$ is
irreducible).

\point Suppose a component of the general curve over $\Xi$ has arithmetic
genus 3 and maps with degree 3.  Then the general curve must have one
other component, with genus 0 and degree 1 (and the two components
meet at one node).  As $\Xi \subset \mts$, the map from a general
curve is a limit of canonical maps.  As the general curve is of
compact type (i.e. the dual graph is a tree), the pullback of
$\oh_{\proj^2}(1)$ to the general curve must be the line bundle
described in the definition of $T$ (see Section \ref{boundarydef}).  Hence
$\Xi \subset T$, so $\Xi = T$ (as $\dim T = 13$ and $T$ is
irreducible).

\begin{figure}
\begin{center}

	   \setlength{\unitlength}{.1\picunit}
	    \begingroup\makeatletter\ifx\SetFigFont\undefined%
\gdef\SetFigFont#1#2#3#4#5{%
  \reset@font\fontsize{#1}{#2pt}%
  \fontfamily{#3}\fontseries{#4}\fontshape{#5}%
  \selectfont}%
\fi\endgroup%
{\renewcommand{\dashlinestretch}{30}
\begin{picture}(7062,3657)(0,-10)
\path(1350,936)(1350,2736)
\path(1275,1986)(2250,3336)
\path(1425,2361)(525,3336)
\path(6150,3336)(6150,936)
\path(7050,3036)(5250,3036)
\put(825,2361){\makebox(0,0)[lb]{\smash{{{\SetFigFont{12}{14.4}{\rmdefault}{\mddefault}{\updefault}$y_2$}}}}}
\put(2100,2736){\makebox(0,0)[lb]{\smash{{{\SetFigFont{12}{14.4}{\rmdefault}{\mddefault}{\updefault}$C_3$ (1,0)}}}}}
\put(0,3486){\makebox(0,0)[lb]{\smash{{{\SetFigFont{12}{14.4}{\rmdefault}{\mddefault}{\updefault}$C_2$ (1,0)}}}}}
\put(1425,1836){\makebox(0,0)[lb]{\smash{{{\SetFigFont{12}{14.4}{\rmdefault}{\mddefault}{\updefault}$y_3$}}}}}
\put(1575,936){\makebox(0,0)[lb]{\smash{{{\SetFigFont{12}{14.4}{\rmdefault}{\mddefault}{\updefault}$C_1$ (2,3)}}}}}
\put(6300,1386){\makebox(0,0)[lb]{\smash{{{\SetFigFont{12}{14.4}{\rmdefault}{\mddefault}{\updefault}$C_1$}}}}}
\put(6300,1011){\makebox(0,0)[lb]{\smash{{{\SetFigFont{12}{14.4}{\rmdefault}{\mddefault}{\updefault}(2,3)}}}}}
\put(4950,3261){\makebox(0,0)[lb]{\smash{{{\SetFigFont{12}{14.4}{\rmdefault}{\mddefault}{\updefault}$C_2$}}}}}
\put(4950,2661){\makebox(0,0)[lb]{\smash{{{\SetFigFont{12}{14.4}{\rmdefault}{\mddefault}{\updefault}(2,0)}}}}}
\put(6225,3111){\makebox(0,0)[lb]{\smash{{{\SetFigFont{12}{14.4}{\rmdefault}{\mddefault}{\updefault}$q$}}}}}
\put(975,36){\makebox(0,0)[lb]{\smash{{{\SetFigFont{12}{14.4}{\rmdefault}{\mddefault}{\updefault}Case i)}}}}}
\put(5700,36){\makebox(0,0)[lb]{\smash{{{\SetFigFont{12}{14.4}{\rmdefault}{\mddefault}{\updefault}Case ii)}}}}}
\end{picture}
} 
\end{center}
\caption{Possibilities for the general map in $\Xi$, Section \ref{g3d2} 
(components labeled (degree, genus))}
\label{Epossibilities}
\end{figure}
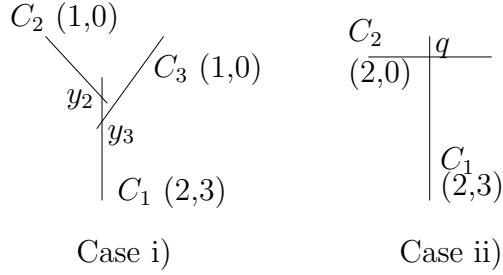

\point
\label{g3d2}
Suppose a component of the general curve has arithmetic genus 3 and
maps with degree 2.  Then the general curve must be one of the
possibilities shown in Figure \ref{Epossibilities}.  In the first
case, $C_1$ meets two components of genus 0, each mapping with degree
1.  In the second case, $C_1$ meets (at one point $q$) a union of
components of total arithmetic genus 0, mapping with total degree 2.

In case i), as the map is a limit of canonical maps, and the source
curve is of compact type, then for some integers $n_2$, $n_3$ the pullback of
$\oh_{\proj^2}(1)$ to $C_1$ is $\cK_{C_1}( (1-n_2)y_2 + (1-n_3) y_3)$,
and the pullback to $C_i$ ($i=2,3$) is $\cK_{C_i}((1+n_i) y_i)$.  From
the degrees of the maps on the components, $n_2=n_3=2$.  As the
pullback of $\oh_{\proj^2}(1)$ to $C_1$ has at least 2 sections,
$$
h^0(C_1, \cK_{C_1}(-y_2-y_3)) \geq 2,
$$
so $\cK_{C_1}(-y_2-y_3)$ must be the
hyperelliptic sheaf, and $y_2$ and $y_3$ must be hyperelliptically
conjugate.  Hence $\Xi \subset Y$, so (as $Y$ is irreducible of
dimension 13) $\Xi=Y$.

In case ii), a similar argument (using $h^0(C_1, \cK_{C_1}(-2q)) \geq 2$)
shows that $q$ is a Weierstrass point of $C_1$.  We claim next that
the image of $C_2$ meets the image of $C_1$ at one point.  Assume
otherwise.  Then the images intersect at two points: the image of $q$,
and some other point $r \in \proj^2$.  (A dimension count shows that
the image of $C_2$ cannot include the image of $C_1$ --- such maps form
a family of dimension less than 13.)  Then consider the germ of this
map above a formal neighborhood of $r$.  The branch of $C_2$ is
immersed in $\proj^2$ and is transverse to the image of $C_1$ (and the two branches
are not connected), so we
can construct the local intersection product of $C_2 \hookrightarrow
\proj^2$ and $C_1
\rightarrow \proj^2$.  These branches intersect with multiplicity 2.
By continuity of intersection products, in any deformation of this
germ of a map the two branches will continue to intersect.  Thus in
any deformation of this germ, the image will remain singular.  Hence
such a map cannot be the limit of smooth maps, so our assumption is
false.  (Remark: this possibility does not appear to be excluded by
the theory of limit linear series.)

Therefore the image of $C_2$ is tangent to the image of $C_1$, so $\Xi
\subset Q$, so (as $Q$ is irreducible of dimension 13) $\Xi = Q$.

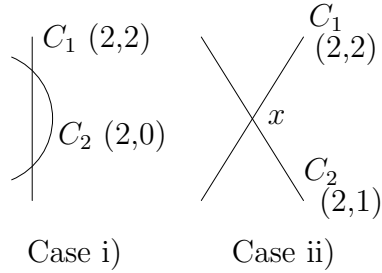
\begin{figure}
\begin{center}

	   \setlength{\unitlength}{.1\picunit}
	    \begingroup\makeatletter\ifx\SetFigFont\undefined%
\gdef\SetFigFont#1#2#3#4#5{%
  \reset@font\fontsize{#1}{#2pt}%
  \fontfamily{#3}\fontseries{#4}\fontshape{#5}%
  \selectfont}%
\fi\endgroup%
{\renewcommand{\dashlinestretch}{30}
\begin{picture}(4892,3657)(0,-10)
\put(-367.000,2136.000){\arc{1950.000}{5.1072}{7.4592}}
\path(308,3336)(308,936)
\path(2783,3336)(4283,936)
\path(2783,936)(4283,3336)
\put(233,36){\makebox(0,0)[lb]{\smash{{{\SetFigFont{12}{14.4}{\rmdefault}{\mddefault}{\updefault}Case i)}}}}}
\put(458,3186){\makebox(0,0)[lb]{\smash{{{\SetFigFont{12}{14.4}{\rmdefault}{\mddefault}{\updefault}$C_1$ (2,2)}}}}}
\put(683,1761){\makebox(0,0)[lb]{\smash{{{\SetFigFont{12}{14.4}{\rmdefault}{\mddefault}{\updefault}$C_2$ (2,0)}}}}}
\put(3233,36){\makebox(0,0)[lb]{\smash{{{\SetFigFont{12}{14.4}{\rmdefault}{\mddefault}{\updefault}Case ii)}}}}}
\put(3758,2061){\makebox(0,0)[lb]{\smash{{{\SetFigFont{12}{14.4}{\rmdefault}{\mddefault}{\updefault}$x$}}}}}
\put(4358,3486){\makebox(0,0)[lb]{\smash{{{\SetFigFont{12}{14.4}{\rmdefault}{\mddefault}{\updefault}$C_1$}}}}}
\put(4283,1236){\makebox(0,0)[lb]{\smash{{{\SetFigFont{12}{14.4}{\rmdefault}{\mddefault}{\updefault}$C_2$}}}}}
\put(4508,786){\makebox(0,0)[lb]{\smash{{{\SetFigFont{12}{14.4}{\rmdefault}{\mddefault}{\updefault}(2,1)}}}}}
\put(4433,3036){\makebox(0,0)[lb]{\smash{{{\SetFigFont{12}{14.4}{\rmdefault}{\mddefault}{\updefault}(2,2)}}}}}
\end{picture}
} 
\end{center}
\caption{Possibilities for the general map in $\Xi$, Section \ref{g2d2} 
(components labeled (degree, genus))}
\label{Fpossibilities}
\end{figure}

\point
\label{g2d2}
Suppose a component of the general curve has arithmetic genus 2 and
maps with degree 2.  Then a quick case check shows that the general
curve must be one of the possibilities shown in Figure
\ref{Fpossibilities}, or the general curve has a contracted union of
components of arithmetic genus 1.  We save the latter case for the end
of the proof, Section \ref{last}.

In case i), $\Xi \subset P$, so $\Xi = P$.  In case F ii), the map
from $C_2$ to $\proj^2$ is given by the line bundle $\cK_{C_2}(2x) \cong
\oh_{C_2}(2x)$, so the double cover from $C_2$ ramifies at $x$.  Hence $\Xi = X$.

We have now completed our list, so we now need to show that there are
no more enumeratively relevant components.

\point
\label{StepG}
Suppose that the general curve has no contracted components, and has
no (arithmetic) genus 2 component mapping with degree 2, and no genus 3 components. 

Replace $\Xi$ by an open subscheme where the topological type of the
source curve is constant.  Then replace $\Xi$ by an \'{e}tale cover where
the components are distinguishable (i.e. the components of the
universal curve correspond to components of a general $k$-fiber).
Let $c$ be the number of irreducible components, and let $\Xi_i$ ($1
\leq i \leq c$) be the families of maps corresponding to the
components of the universal curve over $\Xi$.  It is straightforward to check that $\idim \Xi
\leq \sum_i \idim \Xi_i$.  Let $n$ be the number of nodes connecting
distinct components of the general fiber, and $d_i$ and $g_i$ the
degree and arithmetic genus of the (map from the) $i$th component (so
$\sum_i d_i = 4$, $\sum_i (g_i-1) + n = 2$).  As $\idim \Xi_i \leq 3 d_i  + g_i - 1$,
\begin{eqnarray*}
13 = \idim \Xi & \leq & \sum_i \idim \Xi_i \\
& \leq & 3 \sum_i d_i + \sum_i  (g_i-1) \\
& \leq & 12 + (2-n) \\
&=& 14 - n. 
\end{eqnarray*}
Hence $n=0$ or 1.  If $n=0$, there is only one component, necessarily
of arithmetic genus 3, contradicting the hypothesis of \ref{StepG} that
there are no genus 3 components.  If $n=1$, there are two components.  But then one of the components
must be genus 3, or genus 2 mapping with degree 2, violating the hypotheses of \ref{StepG}.

\point  
\label{clump}
Finally, we show that the general curve of $\Xi$ cannot have any
components contracted by the map.  If $C \rightarrow \proj^2$ is a
stable map, and $D$ is a connected union of contracted components of
$C$ not meeting any other contracted components of $C$, we say $D$ is
a {\em contracted clump}.  Note that if a stable map is ``smoothable''
(i.e. can be deformed to a map from a smooth curve), then any
contracted clump cannot just meet a single, immersed branch --- it
must meet either at least two non-contracted branches, or one
contracted branch $C$ at a point $p$ such that the map $C \rightarrow
\proj^2$ ramifies at $p$.  (More generally, it is also true ---
although not immediate --- that if a stable map is smoothable, a
contracted clump meets the rest of $C$ at one point $p$, the image of
the germ of $C$ at $p$ is reduced, and the map is unibranch over the
image of $p$ (no other branches of $C$ ``interfere'' with the picture)
then the arithmetic genus of the clump is at most the
$\delta$-invariant of the image of the germ of $C$ at $p$.)

\point  Suppose that the general curve of $\Xi$ has at least one  contracted component,
no genus 2 component mapping with degree 2, and no genus 3 components. 

As in \ref{StepG}, reduce to the case where the components of the universal curve over $\Xi$ are
distinguishable.  Let $c$ be the number of components mapping with
positive degree to $\proj^2$.  Base change further if necessary so the nodes of the universal curve over $\Xi$
are also distinguishable.  

Construct the family $\Xi'$ by (i) taking the closure (in the
universal curve over $\Xi$) of the generic points of the
non-contracted components (essentially discarding the contracted
components), and (ii) for every contracted clump meeting more than two
non-contracted branches, choose two of the branches and glue them
together at a node.  (To be precise, the schemes $\Xi$ are $\Xi'$ are
the same, but the families above them are different.)  Then as in 
\ref{clump}, the maps in $\Xi$ through a fixed point (resp. tangent to
a fixed line) are the same as the maps in $\Xi'$ through a fixed point
(resp. tangent to a fixed line).  (The gluing described above was to
ensure that a line tangent to a map in $\Xi$ because it passed through
the image of a contracted clump is also tangent to the corresponding
map in $\Xi'$ because it passes through the image of a node.)  Thus
$\idim(\Xi) = \idim(\Xi')$: the contracted components ``do not
contribute to intersection dimension''.

Next, let $\Xi_i$ ($1 \leq i \leq c$) be the families of maps
corresponding to the components of the universal curve over $\Xi'$, so
$c$ is the number of non-contracted components in $\Xi$.  Let $d_i$
and $g_i$ ($1 \leq i \leq c$) be the degree and genus of the maps in
$\Xi_i$.  Let $b$ be the number of contracted clumps, and $h_1$,
\dots, $h_b$ their arithmetic genera.  Call the nodes on
non-contracted components of the universal curve over $\Xi$ {\em
eligible nodes} (so each eligible node lies on at most one contracted
clump).  Let $n$ be the number of eligible nodes, so $$
\sum_{i=1}^c (g_i - 1) + \sum_{j=1}^b (h_j-1) +n = 2.
$$
Then
\begin{eqnarray*}
13 = \idim \Xi = \idim \Xi' &\leq& \sum_{i=1}^c \idim \Xi_i \\
& \leq & 3 \sum_{i=1}^c d_i + \sum_{i=1}^c (g_i-1) \\
& \leq & 12 + (2-n) + \sum_{j=1}^b (1-h_j),
\end{eqnarray*}
so 
\begin{equation}
\label{contradiction}
n-1 \leq \sum_{j=1}^b (1-h_j).
\end{equation}
For reasons of stability, a contracted clump with arithmetic genus 0 must have at
least 3 eligible nodes.  If $r$ is the number of 
such ``genus 0'' contracted clumps, then the right side of
(\ref{contradiction}) is at most $r$, while the left side is at least
$3r-1$, so $r=0$.

Hence the right side of (\ref{contradiction}) is at
most zero, so $n=0$ or 1.  As (by hypothesis of this step) there is a
contracted component, $n>0$, so $n=1$, and the left side is 0.  Hence
$b=1$ and $h_1=1$, so the map must be from a genus 2 curve $C_1$
mapping with degree 4, union a contracted genus 1 curve $C_2$, meeting
at a single point $p$.  By \ref{clump}, the map from $C_1$ is not an
immersion at $p$.  The intersection dimension of the family is the
same as that of the family of maps from $C_1$ (with the contracted
component discarded), and if this is 13, then from Table \ref{idimmax}
the general map from $C_1$ must be an immersion, giving a
contradiction.

\point
\label{last}
Finally, we take care of the remaining case from \ref{g2d2}, if a
component $C_1$ of the general curve has arithmetic genus 2 and maps
with degree 2, and there is a contracted clump $C_2$ of arithmetic
genus 1 (and at least one more component, for degree reasons).  Then
the genus $g_i=2$ degree $d_i=2$ map moves in a family of
(intersection) dimension at most $(3d_i+g_i-1) + 1$, so the same
argument as in the previous step gives $$ n-2 \leq \sum_{j=1}^b
(1-h_j).  $$ If $r$ is the number of ``genus 0 contracted clumps'',
then the left side is at least $3r-1$ (as there is at least one
eligible node on the genus 1 contracted clump, and at least 3 on each
genus 0), and the right side is at most $r$, so $r=0$.  Hence $b=1$
and $n=2$, and the other non-contracted component must be a rational
curve $C_3$ mapping with degree 2.  The map from $C_1$ moves in a
family of intersection dimension at most 8, and the map from $C_3$
moves in a family of intersection dimension at most 5, so (as $\idim
\Xi = 13$) equality holds in both cases.  For a general $k$-point
in $\Xi$, the image of $C_3$ (a smooth conic) is transverse to the map
of $C_1$ (a line) at two points; let $b_1$ and $b_2$ be these two
(smooth, immersed) branches of $C_3$.  Neither branch can be a smooth
point of the total curve $C_1 \cup C_2 \cup C_3$, as then the map
wouldn't be smoothable by the same argument as \ref{g3d2} Case ii).
Hence one of the branches is a point attached to $C_1$, and the other
is a point of attachment to the collapsed elliptic curve $C_2$ (and
this accounts for both nodes of $C_1 \cup C_2 \cup C_3$).  But then
this contracted clump ($b_2 \cup C_2$) isn't smoothable by Section
\ref{clump}, giving a contradiction.

This completes the proof of Theorem \ref{ugly}.
\epf

\section{Determining coefficients using test families}
\label{testfamilies}

We now determine as many of the unknown co-efficients in (\ref{e1})
and (\ref{e2}) as we can easily, using test families.  (Although they
will not be used here, other methods, such as pencils --- as in
Section \ref{cubicaside} --- and torus actions give test families with
which $h$, $i$, $t$, $p$, $h'$, $i'$, $t'$, $p'$ could be determined.)

Suppose $\pi:  C \rightarrow S$ is a family of nodal genus 3 curves over a one-dimensional smooth base.  For
convenience, let $\om:= \om_{C/S}$.  Let $L$ be
an ample line bundle on $S$.  Suppose $Z$ is a union of  components of fibers, 
and that the total space of $C$ is smooth at all points of $Z$.  Let $M = \om(-Z) \otimes \pi^*L^n$.
Suppose that for every $s \in S(k)$,
\begin{equation}
h^0(\oh_{C_s},\om(-Z)|_{C_s})=3.
\label{telegraph}
\end{equation}
Then $\pi_* M$ is a rank 3 vector bundle on $S$ (Grauert, \cite{hart}
Cor. III.12.9).  Suppose $n  \gg 0$, so $\pi_* M$ is generated by global
sections.  Then for generally chosen sections, all degeneracy loci of
$\pi_* M$ are reduced of the expected dimension (\cite{f} E.g.
14.3.2).  Thus three general sections of $H^0(C, M)$ determines a map
of nodal curves $\rho: C \rightarrow \proj^2 \times S$, and this
linear system is base-point free.

\point Let $I'$ be the (scheme-theoretic) degeneracy locus where the three
sections are linearly dependent; $I'$ is dimension 0, and (as $S$ is
smooth) we will denote the associated (Weil or Cartier) divisor $I'$
as well.  Note that $I'$ and $\pi(Z)$ are disjoint.  Away from the
fibers above $\pi(Z)$ and $I'$, $\rho$ is an immersion.  For the rest
of this section, assume $\rho: C \rightarrow \proj^2 \times S$ is a
family of stable maps whose general curve is smooth (so $C$ has at
worst $A_n$ singularities), inducing a morphism $\phi: S \rightarrow
\mts$.  (A priori the family only induces a morphism $S \rightarrow
\mtp$, but as $S$ is normal, the morphism lifts to
$\phi$.)

Simple calculation using $\phi^* \al = M^2$ and $\phi^* \be = M\cdot
(M+ \om)$ (\cite{ell} 3.10) gives
\begin{eqnarray}
\deg_S \phi^* \al &=& \deg_S(\om^2 - 2 \om Z + Z^2) + 8n \deg_SL,
\label{al} \\
\deg_S \phi^* \be &=& \deg_S(2\om^2 - 3 \om Z + Z^2) + 12n \deg_SL.
\label{be} 
\end{eqnarray}

\tpoint{Proposition}  {\em  If $\pi_* \oh_C(Z) = \oh_S$, and $\eta$ is the locus of nodes of 
the family, then 
$$
\deg_S \phi^* I = 3n \deg_S(L) + \frac 1 {12} \deg_C(\om^2-6\om Z + 6Z^2 + \eta).
$$
}
\label{grr}

Note that $\pi_* \oh_C(Z) = \oh_S$ if $Z$ is a positive linear
combination of components of fibers of $\pi$, and $Z$ does not contain any
fibers of $\pi$.

\bpf By Proposition \ref{iamdegenerate}, $\phi^* I = I'$.  As $I'$ is a degeneracy locus, by \cite{f} Ch. 14,
$$I' = c_1(\pi_*(M)) = c_1(\pi_*  (\om(-Z))) +  3 \deg_S (L^n ).$$
From (\ref{telegraph}) and Serre duality, $h^0(C_s,\oh_C(Z)|_{C_s})=1$ for all $s \in S(k)$, so
$R^1 \pi_* ( \om(-Z)) = (\pi_* \oh_C(Z))^{\vee} = \oh_S$ (by \cite{hm} Ex. 3.12).

By Grothendieck-Riemann-Roch,
\begin{eqnarray*}
ch \pi_* \om(-Z) &= & ch R^1 \pi_* \om(-Z) + \pi_* \left(  \left( 1 - (\om-Z) + \frac {(\om-Z)^2} 2 \right) \cdot \left( 1 - \frac \om 2 +  \frac {\om^2+\eta} {12} \right) \right) \\ 
&=& 3 + \pi_* \left( \frac {\om^2 + 6 \om Z + 6 Z^2 + \eta} {12} \right),
\end{eqnarray*} 
and the result follows after simple manipulation.
\epf

\epoint{Calculating $i$ and $i'$}
Fix a general genus 3 curve $C_1$, and let $C = C_1 \times \proj^1$ and $S = \proj^1$.  Apply the
set-up above with $Z=0$ and $L = \oh_S(1)$.  Then of the divisors appearing in (\ref{e1}) and (\ref{e2}),
only $\al$, $\be$, and $I$ intersect the image of $S$ in $\mts$.  From (\ref{al}), (\ref{be}), and Prop. 
\ref{grr}, $\deg_S \phi^* \al = 8 n$, $\phi^* \be = 12 n$, $\phi^* I = 3n$.  Substituting this
into (\ref{e1}) and (\ref{e2}) (pulled back to $S$) yields $i=12$ and $i'=72$.

\epoint{Calculating $h$ and $h'$} \label{hhprime}
Let $\psi:  S \rightarrow \cmbar_3$ be any morphism from a smooth curve $S$, such that $\psi^* h$ and $\psi^* \de_0$ are non-empty unions 
of reduced points and $\psi^* \de_1$ is empty.  (One such family is given in \cite{hm} Ex. (3.166) part 3.)
 Let $C$ be the pullback of the universal curve to $S$ (so $C$ is smooth).  
Apply the set-up above with $Z =0$ and $L$ any degree 1 (ample) divisor on $S$.  Then the image of $S$ in $\mts$ misses all
divisors in (\ref{e1}) and (\ref{e2}) except $\al$, $\be$, $I$, $H$, and $\De_0$.  From (\ref{al}) and (\ref{be}) and
Prop. \ref{grr}, 
$\deg_S \phi^* \al = \deg_C \om^2 +8n$, $\deg_S \phi^* \be = 2 \deg_C \om^2 + 12 n$, 
and $12 \deg_S \phi^* I = \deg_C \om^2 + 36 n + \deg_S \phi^* \De_0$.  By \cite{hm}
p. 158 and p. 188, 
$$
4h = 3 \pi_* \om^2_{\tC/\cmbar_3} - \de_0 - 9 \de_1
$$
as divisors on the stack $\cmbar_3$ (where $\pi:  \tC \rightarrow \cmbar_3$ is the universal curve).
Now $\psi^* h$ is a union of reduced points, and by Section \ref{m3bar} $\psi^* h = a \phi^* H$, so $a=1$ and $\psi^* h = \phi^* h$.  
Thus 
$$
4 \deg_S \phi^* H = 3 \deg_C \om^2 - \deg_S \phi^* \De_0.
$$ 
Substituting into (\ref{e1}) and (\ref{e2}) yields $h=4$ and $h'=28$.

\epoint{Calculating $t$ and $t'$}

Fix a general genus 3 curve $C_1$, and let $C$ be the blow-up of $C_1 \times \proj^1$ at a general point with
exceptional divisor $E$, and let $S=\proj^1$.  Apply the usual set-up with $Z=2E$ and $L=\oh_S(1)$.
All divisors in (\ref{e1}) and (\ref{e2}) are 0 except $\al$, $\be$, $I$, and $T$.
Then $\deg_S \phi^* \al = -1 + 8n$ and $\deg_S \phi^* \be = 12n$.  Also,
$\deg_S \phi^* T = 1$ by Criterion \ref{node1}.  By Proposition \ref{grr}, 
$12 \deg_S \phi^* I = 36n - 12$.  Substituting into (\ref{e1}) and (\ref{e2}) yields
$t=6$ and $t'=45$.

\epoint{Calculating $p$ and $p'$}
Let $\psi: S \rightarrow \cmbar_3$ be any morphism from a smooth curve
$S$ such that $\psi^* \de_1$ is empty and $\psi^* \de_0$ is a union of
reduced points plus one point $p$ with multiplicity 2.  (For example,
double-cover the base of the family in Section \ref{hhprime} ramified at one of the points
mapping to $\de_0$, and at other generally chosen points.)  Let $C'
\rightarrow S$ be the pullback of the universal curve over $\cmbar_3$,
so $C'$ is smooth except for an $A_1$-singularity above $p$.  Let $b:
C \rightarrow C'$ be the blow-up of $C'$ at the singularity, with
exceptional divisor $E$, so $C$ is smooth and $b^* \om_{C'/S} = \om$.

Apply the usual construction, with $Z=E$, and $L$ a degree 1 (ample) line bundle
on $S$.  The divisors appearing in (\ref{e1}) and (\ref{e2}) intersecting this
family are $\al$, $\be$, $P$, $H$, $I$, and $\De_0$.  One may check that on $C$,
$\om \cdot Z = 0$ and $Z^2 = -2$, so
$\deg_S \phi^* \al = \deg_C \om^2 - 2 + 8n$,
$\deg_S \phi^* \be = 2\deg_C \om^2 - 2 + 12n$,
$\deg_S \phi^* P = 1$ (by Criterion \ref{node1}), and
$12 \deg_S \phi^* I = 36 n + \deg_C \om^2 - 12 + \deg_S \phi^* \De_0 + 2$ (by Proposition \ref{grr}).

From the family $C' \rightarrow S$ (as $\deg_S \psi^* \de_0 = \deg_S \phi^* \De_0 + 2$)
we have (as in Section \ref{hhprime})
\begin{eqnarray*}
4 \deg_S \phi^* H &=& 4 \psi^* h \\ &=& 3 \deg_{C'} \om^2_{C'/S} - (\deg_S \psi^* \de_0+2) \\
&=&3 \deg_C \om^2 - \deg_S \phi^* \De_0 -2
\end{eqnarray*}
Pulling back (\ref{e1}) and (\ref{e2}) to $S$ and solving for $p$ and $p'$ yields $p=2$ and $p'=20$.

\epoint{Calculating $x$ and $x'$}
Let $\psi: S \rightarrow \cmbar_3$ be any morphism from a smooth curve
$S$ such that $\psi^* \de_0$ and $\psi^* \de_1$ are unions of reduced
points and $\psi^* \de_1$ is non-empty.  Let $C$ be the pullback of
the universal curve to $S$ (so $C$ is smooth).  Let $m = \deg_S \psi^*
\de_1$, and let $Z$ be the union of the ($m$) genus 1 components of
fibers.  Apply the usual construction with $L$ a degree 1 (ample) line
bundle on $S$.  All divisors in (\ref{e1}) and (\ref{e2}) are 0 except
$\al$, $\be$, $I$, $\De_0$, $H$, and $X$.  Simple calculations yield
$\deg_C Z^2 = -m$ and $\deg_C \om Z = m$, so $\deg_S \phi^* \al =
\deg_C \om^2 - 3m + 8n$, $\deg_S \phi^* \be = 2 \deg_C \om^2 - 4m + 12
n$.  By Criterion \ref{node1}, $\deg_C \phi^* X = m$ (so $c=1$).  By Proposition
\ref{grr}, $12 \deg_S \phi^*I = \deg_C \om^2 + \deg_S \De_0 -
11m+36n$.  As in Section \ref{hhprime}, $4h = 3 \pi_* \om_{\tilde{C}/\cmbar_3}^2 -
\de_0 - 9 \de_1$ on $\cmbar_3$, so $$ 4 \deg_S \phi^* H = 3 \deg_C
\om^2 - \deg_S \phi^* \De_0 - 9m.  $$ Substituting these values into
(\ref{e1}) and (\ref{e2}) gives $x=6$ and $x'=48$.

\epoint{Aside:  Multiplicities of discriminants}
As a consequence, we can compute the multiplicity of the discriminant
hypersurface $\De$ (in the parameter space $\proj^{14}$ of quartics)
at various points.  Let $p$ be a general point of the locus in
$\proj^{14}$ corresponding to the divisor $H$ (respectively $I$, $T$, $P$).
Then construct a family of maps by taking a general pencil through
$p$.  If $m$ is the multiplicity of the discriminant at $p$, and $a$
is the order of the automorphism group of the limit map, so $a=2$ (resp. 4,3,2) then the pencil
meets $\al$ with degree 1, $\De_0$ with degree $(\deg \De - m) = 27-m$,
and $H$ (resp. $I$, $T$, $P$) with multiplicity $1/a$.  Then from (4),
using $h' = 28$ (resp. $i'=72$, $t'=45$, $p'=20$), the multiplicity of
$\De$ at $p$ is $m=14$ (resp.  18, 15, 10), recovering examples of
Aluffi and F. Cukierman (\cite{ac} Example 3.1).

\section{Characteristic numbers of boundary divisors}
We next calculate the characteristic numbers of the boundary divisors;
the final answers are given in Table \ref{boundary}.
Maple code computing many of the characteristic numbers described here is available upon request.
As Zeuthen had essentially calculated these before (\cite{zeuthen} p. 391,
see Section \ref{zsection}), we have a quick check on our numbers.

\begin{table}
\begin{center}
\begin{scriptsize}
\begin{tabular}{|c|c|c|c|c|c|c|c|c|} \hline
$a$ 	& $\ga_a\De_0$ 	&$\ga_a H$&$\ga_aI$&$\ga_aT$	&$\ga_aP$&$\ga_aQ$	&$\ga_aX$&$\ga_aY$ \\ \hline
13 	&  27		& 0	& 0	& 0		& 0	&  0		& 0	& 0	\\
12	& 162  		& 0	& 0	& 0		& 0	& 0		& 0	& 0	\\
11	& 972		& 0	& 0	& 0		& 0	& 0		& 0	& 0	\\
10	& 5832		& 0	& 0	& 0		& 0	& 0		& 0	& 0	\\
9	& 34992 	& 0	& 0	& 0		& 0	& 0		& 0	& 0	\\
8	& 209952 	& 0	& 0	& 0		& 0	& 0		& 0	& 0	\\
7	& 1256352	& 0	& 0	& 0		& 168	& 0		& 0	& 0	\\
6	& 7453872	& 0	& 0	& 0		& 4536	& 72		& 0	& 0	\\
5	& 43393596	& 4096	&  0	& 0		& 69860	& 1972		& 0	& 150	\\
4	& 242612208	& 110592& 0	& $54 \tau$	&716688	& 24210		& 4032	& 2700	\\
3	&1268876232	&1635840& 0	&$103320 + 1170 \tau$	&5332320&177300		& 105840& 19170	\\
2	&5919651072	&14805120&$16\iota$&$1523720+10120 \tau$	&29220576&842160	& 1164240& 59400	\\
1	&23328812592	&90549360&$288\iota$&$8651280+40920 \tau$&115886232&2561724	&7609140& 0	\\
0	&74651593680	&403572312&$2535\iota$ &0 	&308287980&4487769	&33648615& 0	\\ \hline
\end{tabular}
\end{scriptsize}
\end{center}
\caption{Characteristic numbers of boundary divisors ($\ga_a = \al^a \be^{13-a}$ for brevity)}
\label{boundary}
\end{table}

\bpoint{Everything but $\De_0$}
The characteristic numbers of the components of the families over each
boundary divisor (involving maps of lower genus and/or degree) are
already known.  Then using \cite{ell} Section 3, we can calculate the
characteristic numbers of the boundary divisors.  For the sake of
brevity, we will explicitly calculate one characteristic number for
each boundary divisor, and hope that the general method is clear.

On \cite{zeuthen} p. 390--391, Zeuthen computes the characteristic
numbers of the boundary divisors when $a=2$, $b=11$ as sums, without
further explanation.  Although
his method of computing the summands is different, his summands agree
with the summands computed by this method.  (See Section
\ref{zsection} for a comparison, a glossary of notation, and further
discussion.)  The interested reader can use this method and
use Zeuthen's sums as a check.

\epoint{The divisor $H$}

We count the double covers of conics ramified at 8
points, passing through a fixed general point, and tangent to 12 fixed
general lines.  The double cover is tangent to a line if 
\begin{enumerate}
\item[(i)] a branch point lies on the line, or
\item[(ii)] the image curve is tangent to the line 
(which will give a multiplicity of 2, for the choice of the two branches
to be tangent to the line).
\end{enumerate}
Thus the characteristic number is a sum over non-negative integers $a$, $b$ with $a+b=12$ (where $a$ of the 12 lines are
tangent in the sense of (i) and $b$ are tangent in the sense of (ii)).  Of the $b$ lines, there are $4-a$ pairs such that the conic
passes through the intersection of that pair (thus fixing the conic, up to a finite number of choices).  The double cover branches
at those $4-a$ points, plus at one point of the conic's intersection with each of the remaining $b-2(4-a)$ lines; this accounts for all
$8=(4-a)+(b-2(4-a))$ ramifications.  The number of such maps is a product of several terms:
\begin{itemize}
\item $\binom {12} a$ from the choice of the $a$ lines,
\item $\frac {b!} { (4-a)! 2^{4-a} ( b-2(4-a))!}$ from the choice of the $(4-a)$ pairs of lines,
\item $2^{b-2(4-a)}$ from the choice of intersection of the $b-2(4-a)$ lines with the conic, and
\item the number of conics tangent to $a$ general lines and  through $5-a$ general points (i.e. a characteristic number of plane conics).
\end{itemize}
The multiplicity with which each such map appears is also a product of terms:
\begin{itemize}
\item $\frac 1 2$ from the automorphism of the stable map,
\item $2$ from the choice of pre-image of the fixed point, and
\item $2^a$ from the choice of tangent point to the $a$ lines.
\end{itemize}
Adding these products for $0 \leq a \leq 4$ gives $\al \be^{12} [H] = 90549360$.

\epoint{The divisor $X$}
\label{thedivisorX}

Note that the general map in $X$ has an automorphism group of order 2
(from the genus 1 component).  The divisor on $X$ corresponding to
maps tangent to a line $\ell$ has three components.  The first (resp.
second) is where the genus 1 (resp. genus 2) double cover branches
over $\ell$, but not at a node of the source curve; this divisor
appears with multiplicity 1.  The third divisor is the locus where the node of the
source curve maps to the line, and this divisor appears with
multiplicity 3 (by \cite{ell} Theorem 3.15): two from the node, plus one because the
genus one component ramifies simply over a general line through the node.

We now count the maps in $X$ passing through a fixed
point and tangent to 12 fixed general lines (with appropriate
multiplicities).  There are seven cases to consider.  For convenience,
let $\ell_1$ be the image of the genus 1 component, and $\ell_2$ the
image of the genus 2 component (so $\ell_1$ and $\ell_2$ are lines).

The first case is if none of the 12 lines pass through the image of
the node $\ell_1 \cap \ell_2$, $\ell_1$ passes through the fixed
point, and $\ell_1$ also passes through the intersection $p$ of a pair
of the 12 lines (thus fixing the choice of $\ell_1$).  The line $\ell_2$
passes through the pairwise intersection of two pairs of lines (fixing
$\ell_2$).  The genus 2 cover branches at those two points, plus where
$\ell_2$ intersects 4 other of the 12 lines.  The genus 1 curve
branches at $\ell_1 \cap \ell_2$, the point $p$, and where $\ell_1$
intersects 2 other of the 12 lines.  Note that we have partitioned the
12 lines into 2 (whose intersection is on $\ell_1$), 2 (where the genus
1 cover also branches), $2 \times 2$  (in 2 pairs, whose intersections are on
$\ell_2$, and 4 (where the genus 2 cover also branches).

The degree of this locus is a product of several terms:
\begin{itemize}
\item $\frac 1 2$ from  the automorphism group of the map
\item 2 from the choice of pre-image of the fixed point on the double cover
\item $ \frac 1 2 \binom {12} {2,2,4,2,2}$ from the choice of partition of the 12 lines.
\end{itemize}
Hence this case contributes 623700.

The remaining cases are similar, and are listed in Table \ref{charX}
The $\frac 1 2$ from the automorphisms of the map and the 2 from the
choice of pre-image of the fixed point always cancel, and are omitted
in the table.  The total of the contributions is $\al \be^{12} [X] = 7609140$.

\begin{table}
\begin{center}
\begin{scriptsize}
\begin{tabular}{|c|cccccccc|} \hline
number of lines           & 0 &0 & 2 & 2 & 1 & 1 & 1 & 1 \\
through node              &&&&&&&& \\
\hline
point condition           & 1 & 2 & 1 & 2 & 1 & 1 & 2  & 2 \\
on cover of genus         &&&&&&&& \\
\hline
number of pairwise        &&&&&&&& \\
intersections of lines    & 1 & 2 & 0 & 1 & 1 & 0 & 2 & 1 \\
lines $\ell_1$    &&&&&&&& \\
passes through    &&&&&&&& \\
\hline
number of other &&&&&&&& \\
lines where genus  & 2 & 1 & 3 & 2 & 2 & 3 & 1 & 2 \\
1 cover branches    &&&&&&&& \\
\hline
number of pairwise        &&&&&&&& \\
intersections of    & 2 & 1 & 1 & 0 & 1 & 2 & 0 & 1 \\
lines $\ell_2$    &&&&&&&& \\
passes through    &&&&&&&& \\
\hline
number of other &&&&&&&& \\
lines where genus   & 4 & 5 & 5 & 6 & 5 & 4 & 6 & 5 \\
2 cover branches        &&&&&&&& \\
\hline
multiplicity from      & 1 & 1 & 9 & 9 & 3 & 3 & 3 & 3 \\ 
lines through node              &&&&&&&& \\
\hline
number of parti-    & $\frac 1 2 \binom {12}{2,2,4,2,2}$ & $\frac 1 2 \binom {12} {1,2,2,5,2}$ & $\binom {12} {2,5,3,2}$ &
$\binom {12} {2,2,2,6}$  & $\binom {12} {1,2,2,5,2}$ & $\frac 1 2 \binom {12} {1,3,4,2,2}$ & $\frac 1 2 \binom {12} {1,1,2,2,6}$ &
$\binom {12} {1,2,2,5,2}$ \\
tions of 12 lines & & & & & & & & \\
\hline
total contribution &&&&&&&& \\ 
(product of   & 623700 & 249480 & 1496880 & 748440 & 1496880 & 1247400 & 249480 & 1496880 \\
previous two rows)     & & & & & & & & \\ \hline
\end{tabular}
\end{scriptsize}
\end{center}
\caption{Calculating the characteristic number $\al \be^{12} [X]$ of $X$}
\label{charX}
\end{table}

\epoint{The divisor $Y$}
\label{thedivisorY}

The image of a curve in $Y$ has a point that looks like an
``asterisk''.  The divisor of maps in $Y$ corresponding to maps
tangent to a line $\ell$ includes the locus where the genus 3 curve
branches over $\ell$ (with multiplicity 1), and the locus where the
asterisk lies on $\ell$ (with multiplicity 4:  two from each of the nodes
of the source curve mapping to $\ell$, by \cite{ell} Theorem 3.15).

We count the maps in $Y$ through 2 general points and tangent to 11
general lines.  The two genus 0 components must each pass through one
of the fixed points.  If $m$ is the image of the genus 3 component,
then $m$ must pass through two intersections of pairs of the 11 lines
(and there are $\frac 1 2
\binom {11} {2,2,7} = 1980$ ways of choosing these two pairs).  Of
these two points plus the 7 intersections of $m$ with the remaining
lines, the genus 3 double cover must branch at 8 of them, and the
asterisk must be at the ninth (contributing a multiplicity of 4 or 16,
depending on the number of lines through the asterisk).  Hence the
characteristic number is $\al^2 \be^{11} [Y] = 1980 (2 \times 16 + 7
\times 4) = 59400$.  (A similar calculation appeared in Section
\ref{cubict}.)

\epoint{The divisor $P$}

We count the maps in $P$ tangent to 13 fixed general lines.  For
convenience, let $c$ denote the image of the rational component (a
conic), and $\ell$ the image of the genus 2 component (a line).  Note that
there are two stable maps with the same $c$ and $\ell$ and given
branch points of the double cover of $\ell$ (coming from the choice of which
branch of the cover the conic is glued to).  This will contribute a
factor of 2 to each of our calculations below.

We consider the cases where $x$ of the lines pass through one of the
nodes of the image, and $y$ lines pass through the other ($0 \leq x
\leq y \leq 2$).   Our results are summarized in Tables \ref{charP} and \ref{charP2}.

If $x=y=0$, then the conic $c$ must be tangent to 5 of the lines
(fixing $c$), and the line $\ell$ must pass through 2 intersections of pairs of
lines (fixing $\ell$); the double cover branches at these 2 points, and also where
$\ell$ intersects the 4 remaining lines.  There are $\frac 1 2 \binom
{13} {5,4,2,2}$ ways of partitioning the lines in this way, giving
(along with the factor of 2 described above) a total of 540540.  If
$(x,y)=(0,2)$ or $(2,2)$, the argument is similar.  These three cases
are the first three columns of Table \ref{charP}.

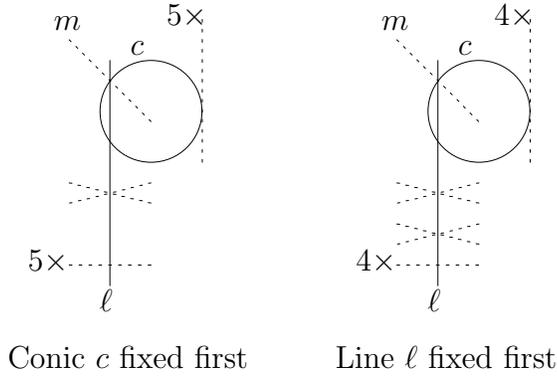
\begin{figure}
\begin{center}

	   \setlength{\unitlength}{.1\picunit}
	    \begingroup\makeatletter\ifx\SetFigFont\undefined%
\gdef\SetFigFont#1#2#3#4#5{%
  \reset@font\fontsize{#1}{#2pt}%
  \fontfamily{#3}\fontseries{#4}\fontshape{#5}%
  \selectfont}%
\fi\endgroup%
{\renewcommand{\dashlinestretch}{30}
\begin{picture}(7924,5205)(0,-10)
\put(2100,3783){\ellipse{1500}{1500}}
\put(6900,3783){\ellipse{1500}{1500}}
\path(1500,4533)(1500,1233)
\dashline{60.000}(900,4833)(2100,3633)
\dashline{60.000}(2850,5133)(2850,3033)
\dashline{60.000}(900,2733)(2100,2433)
\dashline{60.000}(2100,2733)(900,2433)
\dashline{60.000}(900,1533)(2100,1533)
\path(6300,4533)(6300,1233)
\dashline{60.000}(5700,4833)(6900,3633)
\dashline{60.000}(7650,5133)(7650,3033)
\dashline{60.000}(5700,2733)(6900,2433)
\dashline{60.000}(6900,2733)(5700,2433)
\dashline{60.000}(5700,1533)(6900,1533)
\dashline{60.000}(5700,2133)(6900,1833)
\dashline{60.000}(6900,2133)(5700,1833)
\put(675,4983){\makebox(0,0)[lb]{\smash{{{\SetFigFont{12}{14.4}{\rmdefault}{\mddefault}{\updefault}$m$}}}}}
\put(2325,5058){\makebox(0,0)[lb]{\smash{{{\SetFigFont{12}{14.4}{\rmdefault}{\mddefault}{\updefault}$5 \times$}}}}}
\put(1800,4608){\makebox(0,0)[lb]{\smash{{{\SetFigFont{12}{14.4}{\rmdefault}{\mddefault}{\updefault}$c$}}}}}
\put(300,1458){\makebox(0,0)[lb]{\smash{{{\SetFigFont{12}{14.4}{\rmdefault}{\mddefault}{\updefault}$5 \times$}}}}}
\put(1350,858){\makebox(0,0)[lb]{\smash{{{\SetFigFont{12}{14.4}{\rmdefault}{\mddefault}{\updefault}$\ell$}}}}}
\put(5475,4983){\makebox(0,0)[lb]{\smash{{{\SetFigFont{12}{14.4}{\rmdefault}{\mddefault}{\updefault}$m$}}}}}
\put(6600,4608){\makebox(0,0)[lb]{\smash{{{\SetFigFont{12}{14.4}{\rmdefault}{\mddefault}{\updefault}$c$}}}}}
\put(6150,858){\makebox(0,0)[lb]{\smash{{{\SetFigFont{12}{14.4}{\rmdefault}{\mddefault}{\updefault}$\ell$}}}}}
\put(7125,5058){\makebox(0,0)[lb]{\smash{{{\SetFigFont{12}{14.4}{\rmdefault}{\mddefault}{\updefault}$4 \times$}}}}}
\put(5100,1458){\makebox(0,0)[lb]{\smash{{{\SetFigFont{12}{14.4}{\rmdefault}{\mddefault}{\updefault}$4 \times$}}}}}
\put(0,33){\makebox(0,0)[lb]{\smash{{{\SetFigFont{12}{14.4}{\rmdefault}{\mddefault}{\updefault}Conic $c$ fixed first}}}}}
\put(4800,33){\makebox(0,0)[lb]{\smash{{{\SetFigFont{12}{14.4}{\rmdefault}{\mddefault}{\updefault}Line $\ell$ fixed first}}}}}
\end{picture}
} 
\end{center}
\caption{Calculating characteristic numbers of $P$:  the case $(x,y)=(0,1)$}
\label{picturep}
\end{figure}

If $(x,y)=(0,1)$ (so an intersection of $c$ and $\ell$ is required to
lie on some line $m$), there are two possibilities, described pictorially
in Figure \ref{picturep}.  First, $c$ could
be tangent to five of the lines; $\ell$ would pass through one of
the two intersections of $c$ with $m$, and the pairwise intersection of
another pair of the lines; the cover of $\ell$ branches at at the
latter point, and the intersection of $\ell$ with the remaining 5
lines.  There are $\binom {13} {1,5,5,2} = 216216$ ways of
partitioning the lines in such a way, and the other factors involved
are 2 (from the 2 stable maps with the same $\ell$, $c$, and branch
points), 2 (from the choice of intersection of $c$ with $m$), and 2
(the multiplicity from the line $m$ through the node) for a total of
1729728.  In this case, we say that the conic $c$ was {\em fixed first} by
the conditions (and then the choice of $\ell$ was determined using
$c$).

Second, if the line $\ell$ is fixed first, the argument is similar
(see the fifth column of Table \ref{charP} and the second half of Figure
\ref{picturep}).

The case $(x,y) = (1,2)$ breaks into two analogous subcases as well
(first two columns of Table \ref{charP2}).

If $(x,y)=(1,1)$, then the line can be fixed first before choosing the
conic (third column of Table \ref{charP2}), or the conic can be fixed first (fourth column), but
there is one additional case (the last column).  Let $m_1$ and $m_2$
be the two lines such that $c$ and $\ell$ are to intersect once on each line.
The conic $c$ is required to be tangent to 4 of the other lines, and the
line $\ell$ is required to pass through the intersection of 2 others.
(The double cover of $\ell$ is required to branch there, and at the intersection
of $\ell$ with the remaining 5 lines.) The number of ways of
partitioning the lines in this way is $\binom {13} {2,4,2,5}=540540$.
The four tangent lines restrict $c$ to move in a one-parameter family,
and the requirement on $\ell$ restricts $\ell$ to a one-parameter
family.  How often in this (combined) two-parameter family do $c$ and
$\ell$ intersect at two distinct points, one on $m_1$ and one on
$m_2$?  This straightforward enumerative question was addressed (in
much more generality) in \cite{ratell}; we sketch a solution here.  Let
$n_1$ and $n_2$ be two $\proj^1$'s with fixed isomorphisms $n_i \cong
m_i$.  Consider the surface $n_1 \times n_2$.  Let $P$ be the point on
the surface corresponding to the point $m_1 \cap m_2$ in each factor.
As $c$ moves in its one-parameter family, it sweeps out a path in $n_1
\times n_2$ corresponding to pairs of points on $m_1$ and $m_2$; this
path is in class $(4,4)$ --- for a fixed general point on $m_2 \in
\proj^2$, there are 2 conics $c$ tangent to the 4 lines and passing
through the point, and each of those conics intersects $m_1$ in 2
points, and similarly with the roles of $m_1$ and $m_2$ reversed.  The
curve $c$ passes through $P$ with multiplicity 2 (by a similar
argument).  As $l$ moves in a one-parameter family, it sweeps out a
path of pairs of points as well, and this path is in class $(1,1)$,
passing through $P$ with multiplicity 1.  These paths intersect with
multiplicity 8, and it can be checked that the paths intersect at $P$
with multiplicity 2 (corresponding to when both $c$ and $l$ pass
through $m_1 \cap m_2$).  Away from $P$, the two paths intersect at 6
points.  Hence there are 6 configurations where $c$ and $l$ are in the
one-parameter families described above, and intersect at two distinct
points, one on $m_1$ and one on $m_2$.  Thus the factors contributing
in this case are thus 6, 540540 (from
the choice of lines), 2 (the factor described at the beginning of this
note), and 4 (the multiplicity from the two lines passing through the
two nodes), giving a product of 25945920.

The sum of these ten numbers is the characteristic number $\be^{13}[P]
= 308287980$.

\begin{table}
\begin{center}
\begin{scriptsize}
\begin{tabular}{|c|ccccc|} \hline
(x,y)                     & (0,0) & (0,2) & (2,2) & (0,1) & (0,1) \\
\hline
component fixed           & & & &$c$ & $\ell$  \\
first                     &&&&& \\ 
\hline
number of lines           & 5 & 4 & 3 & 5 & 4 \\
tangent to $c$            &&&&& \\ 
\hline
number of pairwise        &&&&& \\
intersections of lines    & 2 & 1 & 0 & 1 & 2 \\
lines $\ell$              &&&&& \\
passes through            &&&&& \\ 
\hline
number of other           &&&&& \\
lines where genus         & 4 & 5 & 6 & 5 & 4 \\
2 cover branches          &&&&& \\
\hline
number of choices         & 1 & 2 & 4 & 1 & 2 \\
for conic $c$             &&&&& \\
\hline
number of choices         & 1 & 1 & 1 & 2 & 1\\
for line $\ell$           &&&&& \\
\hline
multiplicity from         & 1 & 4 & 16& 2 & 2 \\
lines through nodes        &&&&& \\
\hline
number of parti-          & $\frac 1 2 \binom {13} {5,4,2,2}$ & $\binom {13} {2, 4, 2, 5}$ & $\frac 1 2 \binom  {13}{2,2,3,6}$
& $\binom {13} {1,5,5,2}$ & $\frac 1 2 \binom {13} {1,4,4,2,2}$ \\
tions of 13 lines         &&&&& \\
\hline
total contribution        &&&&& \\ 
($2 \times$ product of     & 540540 & 8648640 & 23063040 & 1729728 & 10810800 \\
previous 4 rows)        &&&&& \\ \hline
\end{tabular}
\end{scriptsize}
\end{center}
\caption{Calculating the characteristic number $\be^{13} [P]$ of $P$, part 1}
\label{charP}
\end{table}

\begin{table}
\begin{center}
\begin{scriptsize}
\begin{tabular}{|c|ccccc|} \hline
(x,y)                     & (1,2) & (1,2) & (1,1) & (1,1) & (1,1) \\
\hline
component fixed           & $c$ & $\ell$ & $c$ & $\ell$ & neither \\  
first                     &&&&& \\ 
\hline
number of lines           & 4 & 3 & 5 & 3 & 4 \\
tangent to $c$            &&&&& \\ 
\hline
number of pairwise        &&&&& \\
intersections of lines    & 0 & 1 & 0 & 2 & 1 \\
lines $\ell$              &&&&& \\
passes through            &&&&& \\ 
\hline
number of other           &&&&& \\
lines where genus         & 6 & 5 & 6 & 4 & 5 \\
2 cover branches          &&&&& \\
\hline
number of choices         & 2 & 4 & 1 & 4 & * \\
for conic $c$             &&&&& \\
\hline
number of choices         & 2 & 1 & 4 & 1 & * \\
for line $\ell$           &&&&& \\
\hline
multiplicity from         & 8 & 8 & 4 & 4 & 4 \\ 
lines through nodes        &&&&& \\
\hline
number of parti-          & $\binom {13} {1,2,4,6}$ & $\binom {13} {1,2,3,2,5}$ &
$\binom {13} {2,5,6}$ & $\frac 1 2  \binom {13} {2,2,2,4,3}$ & $\binom {13}{2,4,2,5}$ \\ 
tions of 13 lines         &&&&& \\
\hline
total contribution        &&&&& \\ 
($2 \times$ product of   & 11531520 & 138378240 & 1153152 & 86486400 & 25945920 \\
previous 4 rows)          &&&&& \\ \hline
\end{tabular}
\end{scriptsize}
\end{center}
\caption{Calculating the characteristic number $\be^{13} [P]$ of $P$, part 2}
\label{charP2}
\end{table}

\epoint{The divisor $Q$}
\label{thedivisorQ}

We count the maps in $Q$ tangent to 13 fixed general lines.  This case
is very similar to the case $P$ above.  For convenience, let $c$ denote
the image of the rational component (a conic), and $\ell$ the image of
the genus 3 component (a line).  The divisor on $Q$ corresponding to
maps tangent to a line $m$ has 3 components.  The first is the locus
where the conic $c$ is tangent to $m$.  The second is the locus where
the double cover of $\ell$ branches over $m$, but not at the node of the
source curve.  (Both of these components appear with multiplicity 1.)
The third is the locus where the node of the source curve maps to $m$.
This component appears with multiplicity 3 for the same reason as in
the case $X$ above: two from the node, plus one because the double
cover of $\ell$ ramifies simply over a general line through the image of the
node.

We consider the cases where 0, 1, and 2 lines pass through the image
of the node.  In each of these cases, the conditions can immediately
fix (up to a finite number of choices) one of the two components $c$
or $\ell$ (and then the choice of that component along with the remaining
conditions fix the other component, up to a finite number of choices).
These possibilities are summarized in the first six columns of Table
\ref{charQ}.  (The multiplicity of 1/2 in the last row comes from the
fact that the general map in $Q$ has automorphism group of order 2.)

The one remaining case is the final column of the table.  In this
case, the conic $c$ is tangent to 4 of the lines, restricting $c$ to a
one-dimensional family.  The line $\ell$ passes through the
intersection of a pair of the lines (and the double cover is required
to branch there, as well as where $\ell$ meets 6 more lines),
restricting $\ell$ to a pencil.  The line $\ell$ and the
conic $c$ are required to intersect on the remaining line (call it
$m$), and be tangent there.  We determine how often this happens.

Consider the surface $S = \proj \left( T_{\proj^2} |_m \right)$, i.e. the
$\proj^1$-bundle over $m$ corresponding to the ordered pair (point $p$
on $m$, line through $p$).  This surface is the rational ruled surface
$\eff_1 = \proj( \oh_m \oplus \oh_m(1))$ with Picard group freely
generated by the class $E$ corresponding to ordered pairs (any point
$p$, $m$), and $F$ corresponding to ordered pairs (a fixed point
$p_0$, line through $p_0$), with $E^2 = -1$, $E \cdot F = 1$, $F^2 =
0$.  Define the class $C=E+F$, so $C^2 = 1$; if $p_1$ is a fixed point
of $\proj^2 \setminus m$, $C$ is the class of the set (any point $p$,
line $\overline{p p_1}$).

As $\ell$ moves in a pencil, it describes a curve $B_1$ in $S$
corresponding to (point $\ell \cap m$, $\ell$); this is in class $C$.
As $c$ moves in a one-parameter family, it describes a curve $B_2$ in
$S$ corresponding to (point $p$ on $c \cap m$, tangent line to $c$ at
$p$).  As there are two conics tangent to 4 fixed general lines
through a fixed point $p_0 \in m$, $B_2 \cdot F = 2$.  As there is one
curve tangent to 4 fixed general lines and tangent to $m$, $B_2 \cdot
E = 1$.  Hence $B_2$ is in class $C+2F$, so $B_1 \cdot B_2 = C \cdot
(C+2F) = 3$.  (Of course, one must check that, for general choice of the lines,
all of these intersections are transverse.)

In conclusion, there are 3 ordered pairs $(c,l)$ tangent at a point of
$m$.  Multiplying this by 3 (the multiplicity arising from the line
$m$ passing through the image of the node), $\binom {13} {1,4,2,6}$
(from the ways of partitioning the 13 lines), and $\frac 1 2$ (from
the automorphism group of the general map in $X$), we see that this
case contributes 810810.

\begin{table}
\begin{scriptsize}
\begin{center}
\begin{tabular}{|c|ccccccc|} \hline
number of lines		& 0 & 0 & 2 & 2 & 1 & 1 & 1 \\
through node  &&&&&&& \\
\hline
component fixed           & $c$ & $\ell$ & $c$ & $\ell$ & $c$ & $\ell$ & neither \\  
first                     &&&&&&& \\
\hline
number of lines           & 5 & 4 & 4 & 3 & 5 & 3 & 4 \\
tangent to $c$            &&&&&&& \\ 
\hline
number of pairwise        &&&&&&& \\
intersections of lines    & 1 & 2 & 0 & 1 & 0 & 2 & 1 \\
lines $\ell$              &&&&&&& \\
passes through            &&&&&&& \\ 
\hline
number of other           &&&&&&& \\
lines where genus         & 6 & 5 & 7 & 6 & 7 & 5 & 6 \\
3 cover branches          &&&&&&& \\
\hline
number of choices         & 1 & 1 & 2 & 1 & 1 & 1 & * \\
for conic $c$             &&&&&&& \\
\hline
number of choices         & 2 & 1 & 1 & 1 & 2 & 1 & * \\
for line $\ell$           &&&&&&& \\
\hline
multiplicity from         & 1 & 1 & 9 & 9 & 3 & 3 & 3 \\
lines through node        &&&&&&& \\
\hline
number of parti-          & $\binom {13} {5,2,6}$ & $\frac 1 2 \binom {13} {4,2,2,5}$ & $\binom {13} {2,4,7}$ & $\binom {13} {2,3,2,6}$ &
                          $\binom {13} {1,5,7}$ & $\frac 1 2 \binom {13} {1,3,2,2,5}$ & $\binom {13} {1,4,2,6}$ \\
tions of 13 lines         &&&&&&& \\
\hline
total contribution        &&&&&&& \\ 
($\frac 1 2 \times$ product of     & 36036 & 135135 & 231660 & 1621620 & 30888 & 1621620 & 810810 \\
previous 4 rows)        &&&&&&& \\ \hline
\end{tabular}
\end{center}
\end{scriptsize}
\caption{Calculating the characteristic number $\be^{13} [Q]$ of $Q$}
\label{charQ}
\end{table}

Adding up the seven subtotals gives the characteristic number $\be^{13}[Q] = 4487769$.

\epoint{The divisor $I$}
We count the maps in $I$ tangent to 13 fixed general lines.  Let
$\ell$ be the image of one such map in $I$.  Such maps are in one of
two forms.

The line $\ell$ could pass through the intersection of two pairs of
the 13 lines.  The quadruple cover must ramify at those 2 points, as
well as the points of intersection of $\ell$ with the remaining 9
lines.  (This specifies the canonical quadruple cover, up to a finite number of
choices.)  This number is $\iota$ by definition (see Section \ref{tegoi}).
There are $\frac 1 2 \binom {13} {2,2,9} = 2145$ ways of partitioning
the 13 lines in this case.

On the other hand, the line $\ell$ could pass through the intersection
of a pair of the 13 lines (restricting $\ell$ to a pencil), and
intersect the remaining 11 lines in distinct points; the cover is
required to branch at these 12 points, and be a canonical map.  This
describes a one-parameter family $C$ in $\mbar_{0,12}$ intersecting
$\De_I$ transversely at $\binom {11 } 2 = 55$ points, and missing the
divisors in $S$ (see Section \ref{tegoi} for notation).  Hence the number of
points in $I$ in this family is $C \cdot \left( \frac \iota {11} \De_I
\right) = 5 \iota$.   The number of ways of partitioning the 13 lines 
is $\binom {13} 2$.

Thus the characteristic number $\be^{13} [I]$ is 
$$
\left( 2145 + 5 \binom {13} 2 \right) \iota = 2535 \iota.
$$

\epoint{The divisor $T$}
This case is similar to the previous one.  We count the maps in $T$
tangent to 9 fixed general lines and passing through 4 fixed points.
Let $\ell$ be the image of the genus 3 triple cover, and let $m$ be
the image of the genus 0 component ($\ell$ and $m$ are both lines).
Then $\ell$ must pass through 2 of the 4 points, and $m$ must pass
through the other 2 (so there are 6 ways of partitioning the points
between the components).  The 2 points on $\ell$ contribute a
multiplicity of 3 each (from the 3 possible choices of pre-image of the
point in the triple cover).  The triple cover must branch where $\ell$
meets the 9 lines, and the node of the source curve must map to $\ell
\cap m$, so by Section \ref{tegot} there are $\tau$ points of $T$ satisfying
these conditions.  Thus the characteristic number $\al^4 \be^9 [T]$ is
$6 \times 3^2 \times \tau =54 \tau$.

\bpoint{Characteristic numbers of $\De_0$}

To calculate the characteristic numbers of $\De_0$, we
need to calculate
\begin{enumerate}
\item[$N_a$] the number of degree 4 maps of smooth genus 2 curves through $a$ fixed general points, and
tangent to $13-a$ fixed general lines (the characteristic numbers of
genus 2 quartics),
\item[$N^L_a$] the number of degree 4 maps of smooth genus 2 curves through $a$ fixed general points, and
tangent to $12-a$ fixed general lines, and with the node of the image
lying on another fixed general line, and
\item[$N^p_a$] the number of degree 4 maps of smooth genus 2 curves through $a$ fixed general points, and
tangent to $11-a$ fixed general lines, and with the node of the image
at a fixed general point.
\end{enumerate}
Then, by \cite{ell} Theorem 3.15,
$$
\deg \al^a \be^{13-a} [ \De_0 ] = N_a  + 2 \binom {13-a} 1 N^L_a + 4 \binom {13-a} 2 N_a^p.
$$
(The 2 and 4 come from the multiplicity from the node, and the
binomial coefficients come from the choice of the $13-a$ lines passing
through the node.)  The values of $N_a$, $N^L_a$, and $N^p_a$ appear in \cite{schubert} p. 187 (Section IV).

In \cite{gp}, T. Graber and R. Pandharipande give recursions for the
characteristic numbers of genus 2 plane curves in $\proj^2$, and
computed $N_a$, verifying Zeuthen's degree 4 numbers $N_a$.  Their
method also works for the numbers $N^L_a$ and $N^p_a$ (\cite{tom}, although they have
not explicitly verified Zeuthen's degree 4 numbers for $N_a^L$ and
$N^p_a$).

\section{Linear algebra}  
\label{linearalgebra}
By Section \ref{testfamilies}, equations (\ref{e1}) and (\ref{e2}) can be rewritten
\begin{eqnarray}
6 \al &=& \be + 4H+12I+6T+2P+qQ+6X+yY \label{e1prime} \\
27 \al &=& \De_0+28H+72I+45T+20P+q'Q+48X+y'Y
\label{e2prime}
\end{eqnarray}
(modulo enumeratively irrelevant divisors).  Intersecting these
relations with $\al^a \be^{13-a}$ ($0 \leq a \leq 13$) and using
Table \ref{boundary} yields 28 equations linear in the unknowns $q$,
$q'$, $y$, $y'$, $\iota$, $\tau$, and the characteristic numbers $\deg
\al^a \be^{13-a} [\mts]$.  (Clearly $\deg \al^{14} [\mts] = 1$:  there is one
quartic through 14 general points.)  Solving these equations (with the
aid of Maple) yields $q=6$, $q'=64$, $y=4$, $y'=46$, $\iota = 451440$,
$\tau = 1552$, and the characteristic numbers of smooth quartics:

\tpoint{Theorem} {\em The characteristic number numbers of smooth plane quartics are as given in Table  \ref{char}.}

\begin{table} \begin{center}
\begin{tabular}{|c|c|} \hline
$a$ 	& $\deg \al^a \be^{14-a} [\mts]$ \\ \hline
14	& 1 \\
13 	& 6 \\
12	& 36 \\
11	& 216	\\
10	& 1296 \\
9	& 7776	\\
8	& 46656	\\
7	& 279600 \\
6	&1668096 \\
5	&9840040 \\
4	&56481396 \\
3	&308389896 \\
2	&1530345504 \\
1	&6533946576 \\
0	&23011191144 \\
  \hline
\end{tabular}
\end{center}\caption{Characteristic numbers of smooth plane quartics}
\label{char}
\end{table}

These numbers confirm Zeuthen's predictions (\cite{schubert} p. 187, \cite{zeuthen} p. 391),
and the first ten confirm the calculations of Aluffi (\cite{aquartics}) and van Gastel (\cite{vangastel}).
For unusual consequences of $\iota = 451440$, see \cite{twelvepoints}.

\tpoint{Theorem}  
\label{weilequation}
{\em Modulo enumeratively irrelevant divisors, 
\begin{eqnarray*}
6 \al &=& \be + 4 H + 12 I + 6 T + 2P+6Q+6X+4Y \\
27 \al &=& \De_0 +28H + 72I + 45T + 20P+64Q+48X+46Y.
\end{eqnarray*}}

\section{Comparison with Zeuthen's method}
\label{zsection}

Zeuthen's long article \cite{zeuthen} is devoted to the goal of calculating
the characteristic numbers of smooth plane quartics.  His approach has
many similarities to this one.  Here is a summary based on the
author's understanding of \cite{schubert} p. 184--7 and the french
summary to \cite{zeuthen}, and suggestions by P. Aluffi.  

Zeuthen's aim appears to be to give a general blueprint for all
degrees, and then illustrate it with cubics and quartics.

Note that the dual of a smooth plane quartic has degree 12.  The
parameter space of smooth plane quartics is naturally a dimension 14
locally closed subvariety $\cm$ of $\proj^{14} \times \proj^{90}$
(where the $k$-points correspond to (smooth quartic $C$, dual to $C$);
for Zeuthen $k=\com$ of course).  Let $\cmbar$ be the closure of $\cm$
in $\proj^{14} \times \proj^{90}$.

Not surprisingly, $\cmbar$ has boundary divisors corresponding to the
enumeratively relevant divisors given in Theorem \ref{ugly}.  A
dictionary between our notation and Zeuthen's is given in Table
\ref{dictionary}.

\begin{table} \begin{center}
\begin{tabular}{|l|cccccccccc|} \hline
Notation here &  $\al$     & $\be$  & $\De_0$ & $H$ & $I$ & $T$ & $P$ & $Q$ & $X$ & $Y$
\\ \hline  
Zeuthen's notation & $\mu$ & $\mu'$ & $\al$   & $\vartheta$ & $\nu$ & $\lambda$ & $\xi$ & $\eta$ & $\kappa$ & $\zeta$ 
\\ \hline
\end{tabular}
\end{center}\caption{Notation for analogous divisors on compactifications of $\cm$}
\label{dictionary}
\end{table}

Zeuthen's description of points on the boundary of $\cmbar$ can be 
interpreted in modern language.  For example, a general
point on $\vartheta$ (our $H$) corresponding to a double cover of a
conic branched at eight points is described as twice the class of a
conic with eight ``sommets'' (in \cite{schubert} in German,
``Rangpunkte'') on the conic.  The projection of this point in $\proj^{14}$
is the square of the equation of the conic, and the projection of this
point in $\proj^{90}$ is the square of the equation of the dual conic,
times the equations of the eight lines in the dual plane corresponding
to the lines through the 8 sommets.  In our language, this corresponds
to the fact that lines through the branch points should count for
single tangencies.

Similarly, a general point of $\xi$ (our $P$) is described as having
double sommets at the nodes, corresponding to the fact that lines
through nodes count for two simple tangencies.  A general point of
$\eta$ (our $Q$) has a triple sommet at the singular point (analogous
to the multiplicity of 3 in Section \ref{thedivisorQ}), and the same
is true of $\kappa$ (our $X$, analogous to the multiplicity of 3 in
Section \ref{thedivisorX}).  A general point of $\zeta$ (our $Y$) has
a quadruple sommet at the singular point (analogous to the
multiplicity of 4 in Section \ref{thedivisorY}).

Using these multiplicities, Zeuthen appears to calculate the characteristic numbers of the boundary
divisors in the same way as described here.  However, rather than
having two unknowns $\iota$ and $\tau$ (from the divisors $I$ and
$T$), he has five unknowns $H=\tau$, $I=3280$, $K=5 \tau + 1640$,
$L=\iota$, and $M = 6 \iota$ corresponding to various enumerative
problems (from the analogous divisors $\nu$ and $\lambda$).  

Zeuthen gives equations analogous to Theorem \ref{weilequation} (\cite{zeuthen}, p. 389):
\begin{eqnarray}
\mu' &=& 6 \mu - 2 \xi - 3 \eta - 4 \zeta - 3 \kappa - 6 \lambda - 12 \nu - 2 \vartheta \label{z1}
\\
\al &=& 27 \mu - 20 \xi - 32 \eta - 46 \zeta - 24 \kappa - 45 \lambda - 72 \nu - 14 \vartheta.
\label{z2}
\end{eqnarray}
The coefficients of $\eta$, $\kappa$, and $\vartheta$ (corresponding
to $Q$, $X$, and $H$) are half the analogous coefficients in Theorem
\ref{weilequation}, because the isotropy group of the generic point of
those divisors on $\mts$ is $\zed/2$ (the general such map has an
automorphism group of order 2).  Zeuthen's characteristic numbers differ from those in Table
\ref{boundary} for the same reason.  (Thus equations (\ref{z1}) and
(\ref{z2}) can be interpreted as equality on the coarse moduli scheme
of $\mts$, modulo enumeratively irrelevant divisors.)

It is not clear to the authour how Zeuthen obtained the co-efficients in
(\ref{z1}) and (\ref{z2}), which is the crux of the calculation.  He certainly does not 
provide details of what he considered routine calculations.  P. Aluffi has pointed out the following
intriguing passage (\cite{zeuthen} p. XI Section 26 of the summary):
\begin{verse}
Ayant trouv\'{e} ... les ordres des distances et des angles infiniment
petits qui s\'{e}parerent les points et les tangents des courbes
singuli\`{e}res de ceux de leurs courbes voisines, nous pouvons faire
usage de cette r\`{e}gle pour d\'{e}terminer directement les
coefficients des formules...
\end{verse}
Aluffi suggests that he may have determined co-efficients by computing
angles (or orders of vanishing of angles), and the detailed figures
at the end of the article seem to corroborate this.

\point 
\label{weirdIfact}
There is one small (but interesting) point where Zeuthen is not
correct (without throwing off his calculation).  One of his unknowns
corresponds to the number of solutions to the following problem
(\cite{zeuthen} p. XXII): given a line in the plane and 11 general
sommets on the line, how many choices of a twelfth sommet are there so
the resulting configuration lies in $\nu \subset \proj^{14} \times
\proj^{90}$ 
(corresponding to our $I$)?  In more modern language, given 11 points
on a line, how many choices are there for a twelfth so there is a
canonical cover (of genus 3, degree 4) branched at those 12 points?
This is a (slightly) different question from that asked in Section
\ref{tegoi}: given 11 general points on a line, how many canonical
covers are there branched at those 11 points?

In fact, for each general genus 3, degree 4 canonical cover, there
are 119 other canonical maps branched at the same points!  In other
words, the natural rational map from the (11-dimensional) space of
smooth genus 3 canonical covers of a line $L$, to its image in
$\Sym^{12} L \cong \proj^{12}$ is not birational, as one would naively
expect, but of degree 120!  (This is because the corresponding divisor
in $\Sym^{12} L$ is very unusual.  This divisor will be discussed
further in \cite{twelvepoints}.)  Zeuthen gives the answer to his question as 451440 points.  The actual
answer is 3762 points, each with multiplicity 120.

\end{document}